\documentclass[a4paper,11pt,twoside]{amsart}

\usepackage{amsmath}
\usepackage{amsthm}
\usepackage{amssymb}
\usepackage[all]{xy}

\newtheorem{thm}{Theorem}[section]
\newtheorem{prop}[thm]{Proposition}
\newtheorem{lem}[thm]{Lemma}
\newtheorem{cor}[thm]{Corollary}
\newtheorem{conjecture}[thm]{Conjecture}

\numberwithin{equation}{section}

\theoremstyle{definition}
\newtheorem{definition}[thm]{Definition}
\newtheorem{remark}[thm]{Remark}
\newtheorem{ex}[thm]{Example}

\newcommand{\im}{\operatorname{im}}

\newcommand{\OO}{{\rm O}}
\newcommand{\Db}{{\rm D}^{\rm b}}
\newcommand{\D}{{\rm D}}
\newcommand{\C}{{\rm C}}
\newcommand{\Km}{{\rm Km}}
\newcommand{\Eq}{{\rm Eq}}
\newcommand{\Aut}{{\rm Aut}}

\newcommand{\NS}{{\rm NS}}
\newcommand{\Pic}{{\rm Pic}}
\newcommand{\ch}{{\rm ch}}

\newcommand{\rk}{{\rm rk}\,}
\newcommand{\coh}{{\cat{Coh}}}

\newcommand{\Mod}[1]{{\ko_{#1}\text{-}\cat{Mod}}}

\newcommand{\Hom}{{\rm Hom}}

\newcommand{\td}{{\rm td}}
\newcommand{\ob}{{\rm Ob}}
\newcommand{\Stab}{{\rm Stab}}
\renewcommand{\dim}{{\rm dim}\,}

\newcommand{\iso}{\cong}
\newcommand{\lto}{\longrightarrow}

\newcommand{\id}{{\rm id}}

\newcommand{\rest}[1]{|_{#1}}
\newcommand{\dual}{^{\vee}}

\newcommand{\epi}{\twoheadrightarrow}
\newcommand{\mor}[1][]{\xrightarrow{#1}}
\newcommand{\isomor}{\mor[\sim]}
\newcommand{\abs}[1]{\lvert#1\rvert}

\newcommand{\cat}[1]{\begin{bf}#1\end{bf}}
\newcommand{\Ext}{{\rm Ext}}
\newcommand{\FM}[1]{\Phi_{#1}}

\newcommand{\R}{\mathbf{R}}
\newcommand{\lotimes}{\stackrel{\mathbf{L}}{\otimes}}

\newcommand{\K}{{\rm Kom}}

\newcommand{\cal}{\mathcal}
\newcommand{\ka}{{\cal A}}

\newcommand{\kc}{{\cal C}}

\newcommand{\ke}{{\cal E}}
\newcommand{\kf}{{\cal F}}
\newcommand{\kg}{{\cal G}}
\newcommand{\kh}{{\cal H}}
\newcommand{\kk}{{\cal K}}

\newcommand{\kn}{{\cal N}}
\newcommand{\km}{{\cal M}}
\newcommand{\ko}{{\cal O}}
\newcommand{\kp}{{\cal P}}

\newcommand{\ks}{{\cal S}}
\newcommand{\kt}{{\cal T}}

\newcommand{\NN}{\mathbb{N}}
\newcommand{\ZZ}{\mathbb{Z}}
\newcommand{\QQ}{\mathbb{Q}}
\newcommand{\RR}{\mathbb{R}}
\newcommand{\CC}{\mathbb{C}}

\newcommand{\HH}{\mathbb{H}}
\newcommand{\PP}{\mathbb{P}}

\newcommand{\grf}[1]{\mbox{$\left \{ #1 \right \}$}}        

\def\cc{\chi}

\def\ee{\varepsilon}

\def\jj{\varphi}

\begin{document}

\title{Automorphisms and autoequivalences of generic analytic K3 surfaces}

\author{Emanuele Macr\`{i} and Paolo Stellari}

\address{E.M.: Hausdorff Center for Mathematics, Mathematisches
Institut, Universit{\"a}t Bonn, Beringstr.\ 1, 53115 Bonn,
Germany} \email{macri@math.uni-bonn.de}

\address{P.S.: Dipartimento di Matematica ``F. Enriques'',
Universit{\`a} degli Studi di Milano, Via Cesare Saldini 50, 20133
Milano, Italy} \email{Paolo.Stellari@mat.unimi.it}



\begin{abstract} This is a systematic exposition
of recent results which completely describe the group of
automorphisms and the group of autoequivalences of generic analytic
K3 surfaces. These groups, hard to determine in the algebraic case,
admit a good description for generic analytic K3 surfaces, and are in fact seen to be closely interrelated.\end{abstract}

\maketitle

\section{Introduction}\label{sec:intro}

In this survey paper we present some recent results about the group
of automorphisms and the group of autoequivalences of generic
analytic K3 surfaces. Recall that a surface of this type is a
non-algebraic K3 surface with trivial Picard group. They are dense
in the moduli space of K3 surfaces and their geometry is
considerably simpler when compared to that of algebraic K3
surfaces. This is mainly due to the lack of curves.

The geometry of K3 surfaces has been widely studied
during the last forty years with particular attention to the
algebraic case. Deep results such as the Torelli Theorem have shown that many interesting geometric properties of these surfaces can be
determined by lattice theory and Hodge structures. For example, the
cup product and the weight-two Hodge structure defined on the second
cohomology group of a K3 surface $X$ determine it uniquely up to
isomorphisms. Thus it forms a key tool in the study of the
automorphism groups of K3 surfaces.

The algebraic case happens to be extremely rich and, in fact, we
do not have a complete description of the automorphism group of algebraic
K3 surfaces at this point. It turns out that the isomorphism class of the Picard
lattice can deeply influence this group. One can easily produce
examples where this group is trivial, finite, infinite and cyclic or
infinite but not cyclic (Section \ref{subsec:alg}). Very often the
geometry of the surface and not just the lattice and Hodge
structures of its second cohomology group gives significant clues to
approach a description of this group. From this intricate situation,
a quite large number of interesting results emerged in the algebraic
setting.

The case of generic analytic K3 surfaces is very different. In
particular the description of the automorphism group is very neat
and only two possibilities occur: either the group is trivial
or the group is cyclic and infinite (see Theorem \ref{thm:sumnonalg}
and \cite{Og1}). The latter case is quite rare and the generic
analytic K3 surfaces where this is realized have been
completely described by Gross and McMullen (see Section
\ref{subsec:McMul} and \cite{GMc,Mc}).

There has been a heightened interest in K3 surfaces recently because of advances in string theory.
Indeed, these surfaces are examples of
Calabi--Yau manifolds and irreducible holomorphic symplectic
manifolds at the same time and their geometry is relatively well
understood. Hence they have become one of the favorite
examples to test the homological interpretation of Mirror Symmetry
due to Kontsevich (\cite{K} but see also \cite{Do2,Hu1,Sd}).

Roughly speaking, this conjecture predicts that, for two Calabi--Yau manifolds
$X$ and $Y$ which are mirror, the bounded derived category of
coherent sheaves on one variety should be equivalent to the derived Fukaya
category of the other one. The bounded derived category of coherent sheaves should
take into account the complex structure of a variety while the
derived Fukaya category should encode the symplectic geometry. Indeed the
objects of the latter category are special Lagrangian submanifolds
with a suitable notion of morphisms obtained by means of Floer
homology.

It is known that the abelian category of coherent sheaves $\coh(X)$
on a smooth projective variety $X$ determines the variety uniquely
up to isomorphism. This is no longer true if we consider the bounded
derived category of $\coh(X)$. More precisely, Bondal and Orlov
proved in \cite{BO1} that if either the canonical or the
anticanonical bundle of $X$ is ample then $X$ is determined uniquely
up to isomorphism by the bounded derived category of $\coh(X)$. But
as soon as the canonical bundle becomes trivial (as for K3 surfaces
and, more generally, for Calabi--Yau manifolds), a large number
of interesting phenomena can happen. Just restricting to the case of algebraic
K3 surfaces, one can produce examples of K3 surfaces $X$ and
arbitrarily many non-isomorphic K3 surfaces $Y$ with the same
derived category (\cite{Og,St}). Nevertheless the number of
isomorphism classes of K3 surfaces with such a close relationship to
$X$ is always finite (\cite{BM}). Conjecturally the number of
isomorphism classes of smooth algebraic varieties with the same
derived categories is always expected to be finite.

One can wonder if the same results are still true when $X$ is
non-algebraic but we will see that this is not the case. As a first
step, we will observe that, in order to be sure that all functors and
dualities we need are well-defined, we should work with a suitable
subcategory of the derived category of $\ko_X$-modules instead of
the bounded derived category of coherent sheaves (Section
\ref{subsec:derivcomplexmnflds }). The main problem in this sense is
the validity of Serre duality. Notice that in our geometric
setting, just involving complex non-projective surfaces, this is not
a stumbling block since, exactly as in the algebraic case, this special
subcategory is naturally equivalent to the bounded derived category
of coherent sheaves (see Proposition \ref{prop:linkcoh} and
\cite{BB}) although some problems arise when we deal with the
derived category of the product of two K3 surfaces. Moreover, the
geometric similarities with the algebraic setting are evident in many
other useful instances. For example, any coherent sheaf admits a
locally free resolution of finite length. Hence the Grothendieck
group of these surfaces is generated by the classes of locally free
sheaves, we can easily define Chern classes for coherent sheaves and
the Riemann-Roch Theorem in its non-algebraic version holds true for
any coherent sheaf. More generally we will recall how
Grothendieck-Riemann-Roch Theorem was proved in \cite{OTT2} after
Chern classes were introduced in \cite{OTT1} for coherent sheaves on
analytic varieties. A somewhat sketchy definition of these Chern classes will be included in
Section \ref{subsec:coh}.

Assuming the point of view of physics one can see that many
interesting dualities in string theory can be reinterpreted as
equivalences (or autoequivalences) of derived categories. This
motivates the growing interest in the study of the group of
autoequivalences of the derived categories of K3 surfaces or, more
generally, of varieties with trivial canonical bundle.

While a very nice and complete description of the group of
autoequivalences is now available for abelian varieties (\cite{Or3})
the same is not true for algebraic K3 surfaces. Even for very
explicit examples where the geometry of the surface is very well
understood, this group remains mysterious.

A conjectural description of it was first proposed by
Szendr\H{o}i in \cite{Sz}. Indeed, Orlov (\cite{Or}) proved, using
results of Mukai (\cite{Mu2}), that any autoequivalence of the
derived category of a smooth projective K3 surface induces an
isomorphism of the total cohomology groups preserving some special
lattice and Hodge structures (Section \ref{subsec:actcohom}).
Szendr\H{o}i observed that such an isomorphism should preserve the
orientation of a particular four-dimensional subspace of the real
cohomology (Conjecture \ref{conj:orient}). So we would expect a
surjective morphism from the group of autoequivalences onto a group of
orientation preserving isometries of the total cohomology group of
an algebraic K3 surface. Notice that this is in perfect analogy with
the Strong Torelli Theorem where the preservation of the K\"{a}hler
cone is required.

On the other hand, until recently no clue about the kernel of this
morphism was available. The real breakthrough in this direction is
due to Bridgeland and to his notion of stability condition on
triangulated categories (\cite{B2}). Roughly speaking a stability
condition is determined by a $t$-structure on the triangulated
category and a $\CC$-linear function on the heart of this
$t$-structure satisfying special properties which imitate the usual
definition of $\mu$-stability for coherent sheaves. Bridgeland's
motivation was to axiomatize the notion of $\Pi$-stability which emerged
from Douglas' work (\cite{Dou}). In particular, the manifold
parametrizing stability conditions on the derived category $\cat{D}$
of some Calabi--Yau manifold quotiented by the natural action of the
autoequivalence group of $\cat{D}$ models an extended version of the stringy K\"{a}hler
moduli space of the Calabi--Yau manifold.

Studying stability conditions on the bounded derived category of
coherent sheaves on smooth projective K3 surfaces, Bridgeland
conjectured that the kernel of the morphism considered by
Szendr\H{o}i should be the fundamental group of some period domain
associated to the variety parametrizing stability conditions on the
derived category (\cite{B}).

A first attempt to verify this conjectural description of the
autoequivalence group due to Bridgeland and Szendr\H{o}i was made
in \cite{HMS} where this group was described for generic twisted
K3 surfaces (not considered in this paper) and generic analytic K3
surfaces. Very surprisingly, the groups of automorphisms and
autoequivalences which are quite mysterious in the algebraic case,
can be explicitly studied in the generic analytic case. Moreover
it turns out that these groups are very closely related in the
non-algebraic setting since the group of autoequivalences of
Fourier--Mukai type of a generic analytic K3 surface just consists
of shifts, twists by a special object called spherical and
automorphisms (Theorem \ref{prop:autgennonproj}). (Notice that the
techniques used in \cite{HMS} to prove this result, and briefly
recalled in the second part of the present paper, are applied in
\cite{HMS1} to prove Szendr\H{o}i's conjecture in the smooth
projective case.)

The aim of this paper is to describe these two groups in full detail,
collecting results contained in many papers and trying to give a
comprehensive, systematic and integrated approach.

\medskip

The paper is organized as follows. Section \ref{sec:geom} deals with
the geometry of K3 surfaces. We briefly recall some results in
lattice theory and give some examples of algebraic and
non-algebraic K3 surfaces. Particular attention is paid to the
examples of non-algebraic K3 surfaces proposed by McMullen. We also
discuss the lattice and Hodge structures defined on the second and
total cohomology groups of K3 surfaces.

Section \ref{sec:automorphisms} is devoted to the complete
description of the automorphism group of generic analytic K3
surfaces. In particular we present the works of Gross, McMullen and
Oguiso. These results are compared through examples with the case of
algebraic K3 surfaces.

In Section \ref{sec:cohder} we identify the appropriate triangulated
category where all functors with geometrical meaning can be derived
(as a special case we consider the Serre functor). We recall a
result of Verbitsky about the abelian category of coherent sheaves
on generic analytic K3 surfaces and state that any coherent sheaf
admits a locally free resolution of finite length. Finally, Chern classes for coherent sheaves
on smooth complex varieties are defined according to \cite{OTT2}.

In the last section we prove that a smooth complex variety with
the same derived category as a generic analytic K3 surface is a
surface of the same type and state Szendr\H{o}i's conjecture for
any (projective or not) K3 surface. We discuss the notion of
stability condition due to Bridgeland and use it to fully describe
the group of autoequivalences of generic analytic K3 surfaces.

\medskip

A few warnings are in order at this point. First of all, there are no
new results in this presentation and we have deliberately omitted many
nice results concerning algebraic K3 surfaces. For example, the
section about automorphisms of algebraic K3 surfaces is mainly based
on simple examples and none of the beautiful deep results in
the literature are considered. This is justified since the algebraic setting is out of the scope of this paper. Secondly,
the proofs are described in detail only when the length and the
complexity are reasonable. When we feel that it may be worthwhile to give the reader an impression of the techniques used in a given proof which is otherwise too involved, we simply outline the main steps. The
final warning concerns the fact that the reader is supposed to have
some familiarity with the geometry of K3 surfaces and with the
language of derived categories. Although some introductory sections
are present in this paper, the less prepared reader is encouraged to
have a look at \cite{Ge}, \cite{BPV} and at the first few chapters of
\cite{H}.

\section{The geometry of K3 surfaces}\label{sec:geom}

In this section we discuss a few basic geometric properties of algebraic
and generic analytic K3 surfaces. After describing some simple
examples, we will introduce the Hodge and lattice structures defined
on the second and on the total cohomology groups of any K3 surface.
The high points of the theory are certainly the Torelli Theorem and
the theorem on the surjectivity of the period map.

While it is very easy to describe explicit examples of algebraic K3
surfaces, the non-algebraic ones are usually presented with a rather
abstract approach (e.g.\ invoking the surjectivity of the period
map). At the end of this section, we present beautiful examples of
complex non-algebraic surfaces that were first described by McMullen.

For sake of completeness we include a short section where we
collect well-known results in lattice theory.

\bigskip

\subsection{Lattices}\label{subsec:lattices}~

\medskip

Lattice theory plays a special role in the study of the geometry of
K3 surfaces. As we will see in a while, interesting aspects of
the geometry of K3 surfaces are encoded by the Hodge and lattice
structures defined on their cohomology groups. Due to a celebrated
result of Orlov (see Theorem \ref{thm:coh}) the lattice theoretic
properties of the cohomology of a K3 surface are also important to
detect equivalences between derived categories.

Therefore, not surprisingly, in this paper we will often need a few
algebraic results in lattice theory. We list them in this section
for the convenience of the reader. The proofs can be found in
\cite{Ni} (see also \cite{Do1} for a readable account).

First of all, recall that a \emph{lattice} is a free abelian group
$L$ of finite rank with a non-degenerate symmetric bilinear form
$b:L \times L\to\ZZ$. For simplicity, we always write $x\cdot
y:=b(x,y)$ and we denote by $\OO(L)$ the group of isometries of a
lattice $L$ (i.e.\ the automorphisms of $L$ preserving the
bilinear form).

\begin{ex}\label{ex:lattices} (i) The \emph{hyperbolic
lattice} $U$ is the free abelian group $\ZZ\oplus\ZZ$ endowed with
the bilinear form represented by the matrix
\[
\left(\begin{array}{cc} 0 & 1 \\ 1 & 0 \\\end{array}\right),
\]
with respect to the standard basis.

(ii) Consider the free abelian group $\ZZ^{\oplus 8}$ with the
bilinear form
\begin{equation*}
\Small{\left( \begin{array}{rrrrrrrr}
2 & 0 & -1 & 0 & 0 & 0 & 0 & 0 \\
0 & 2 & 0 & -1 & 0 & 0 & 0 & 0 \\
-1 & 0 & 2 & -1 & 0 & 0 & 0 & 0 \\
0 & -1 & -1 & 2 & -1 & 0 & 0 & 0 \\
0 & 0 & 0 & -1 & 2 & -1 & 0 & 0 \\
0 & 0 & 0 & 0 & -1 & 2 & -1 & 0 \\
0 & 0 & 0 & 0 & 0 & -1 & 2 & -1 \\
0 & 0 & 0 & 0 & 0 & 0 & -1 & 2 \\
\end{array} \right).}
\end{equation*}
Such a lattice is denoted by $E_8$.

(iii) We denote by $\langle d\rangle$ the lattice $L=(\ZZ,b)$ where
$b(h,h)=d$ for a generator $h$ of the underlying rank one abelian
group.

(iv) Given a lattice $L$ with bilinear form $b_L$ and an integer
$m\in\ZZ$, the lattice $L(m)$ coincides with $L$ as a group but its
bilinear form $b_{L(m)}$ is such that
$b_{L(m)}(l_1,l_2)=mb_L(l_1,l_2)$, for any $l_1,l_2\in L$.\end{ex}

A lattice $(L,b)$ is \emph{even} if, for all $x\in L$,
$x^2:=b(x,x)\in 2 \mathbb{Z}$, while it is \emph{odd} if there
exists $x\in L$ such that $b(x,x)\not \in 2 \mathbb{Z}$. Given an
integral basis for $L$, we associate to the bilinear form $b_L$ a
symmetric $(\rk L\times\rk L)$-matrix $S_L$. The {\it discriminant}
of $L$ is the integer $\mathrm{disc}(L):=\mathrm{det} \, S_L$. A lattice
is \emph{unimodular} if $\mathrm{disc}(L)=\pm1$ (notice that $\mathrm{disc}(L)\neq 0$, because $L$ is non-degenerate). Moreover, a
lattice $L$ is positive definite if $x^2>0$, for any $0\neq x\in L$
while $L$ is negative definite if $x^2<0$, for any $0\neq x\in L$.
It is indefinite if it is neither positive
definite nor negative definite. (Notice that $L(-1)$ is negative
definite if and only if $L$ is positive definite.)

Unimodular indefinite lattices are particularly relevant for the
study of K3 surfaces and they are completely classified by the
following result.

\begin{thm} {\bf (\cite{Mi})}\label{thm:milnor} Let $L$ be a unimodular indefinite
lattice. If $L$ is odd then $L\cong\langle 1\rangle^{\oplus m}\oplus
\langle -1\rangle^{\oplus n}$, while if $L$ is even $L\cong
U^{\oplus m}\oplus E_8(\pm 1)^{\oplus n}$, for a unique pair of positive
integers $m$ and $n$.\end{thm}

We can associate a finite group to $L$. Indeed, consider
\[
L^\vee:=\Hom_\ZZ(L,\ZZ)\cong\{l\in L\otimes\QQ:l\cdot p\in\ZZ,\mbox{
for any }p\in L\}
\]
and the natural inclusion $L\hookrightarrow
L^\vee$. The bilinear form on $L$ induces a symmetric bilinear form $b^\vee: L^\vee \times L^\vee \to \QQ$. The \emph{discriminant group of $L$} is the quotient
$$A_L:=L^\vee/L.$$ The order of $A_L$ is $|\mathrm{disc}(L)|$ (see \cite[Lemma
2.1]{BPV}). Moreover, $b^\vee$ induces a
symmetric bilinear form $b_L:A_L \times A_L \rightarrow
\mathbb{Q}/\mathbb{Z}$ and thus a quadratic form
$$q_L:A_L\rightarrow\mathbb{Q}/\mathbb{Z}.$$ Whenever $L$ is even, $q_L$ takes values in
$\QQ/2\ZZ$. By $\OO(A_L)$ we denote the group of automorphisms of
$A_L$ preserving $q_L$. We denote by $\ell(L)$ the minimal number of
generators of $A_L$. The inclusion of $L$ into $L^\vee$ yields a
homomorphism $\Psi:\OO(L)\rightarrow\OO(A_L)$.

The discriminant group of a lattice $L$ and the group of its
automorphisms preserving the quadratic form just defined always play
a role when we want to extend isometries of sublattices to isometries
of $L$ itself. Since at some point we will have to do this, we
briefly recall the main results proved by Nikulin in this direction.

\begin{prop}\label{prop:surj} {\bf (\cite{Ni}, Theorem 1.14.2.)} Let $L$ be an even
indefinite lattice such that $\rk L\geq\ell(L)+2$. Then the map
$\Psi:\OO(L)\rightarrow\OO(A_L)$ is onto.\end{prop}

An embedding $i:V\hookrightarrow L$ of a lattice $V$ into a lattice
$L$ is \emph{primitive} if $L/i(V)$ is a free abelian group. One can
show that given a primitive sublattice $K$ of a unimodular lattice
$L$, there exists a morphism of groups $\gamma:A_K\longrightarrow
A_{K^\perp}$. Moreover,

\begin{prop}\label{prop:ort} {\bf (\cite{Ni}, Corollary 1.6.2.)} If $K$ and $L$ are as above, then the morphism $\gamma$
induces an isometry
$\gamma:(A_K,q_K)\rightarrow(A_{K^\perp},-q_{K^\perp})$.\end{prop}

Hence, $\gamma$ defines an isomorphism
$\psi:\OO(A_K)\rightarrow\OO(A_{K^\perp})$. Furthermore, given an
even unimodular lattice $L$, there are two homomorphisms
$\Psi_1:\OO(K)\rightarrow\OO(A_K)$ and
$\Psi_2:\OO(K^\perp)\rightarrow\OO(A_{K^\perp})$ naturally
associated to a primitive sublattice $K$ of $L$.

\begin{prop}\label{prop:ext} {\bf (\cite{Ni}, Theorem 1.6.1 and Corollary
1.5.2.)} An isometry $f\in\OO(K)$ lifts to an isometry in $\OO(L)$
if and only if $\psi(\Psi_1(f))\in\im(\Psi_2)$.\end{prop}

For later use we include here a simple example:

\begin{ex}\label{ex:latt1} Consider a primitive sublattice
$L=\langle2d\rangle$ of $\Lambda$, with $d>0$. From the
definition we have that $A_L\cong\ZZ/2d\ZZ$. If $d \neq 1$, then it was proved in \cite{Sca} that
$\OO(A_L)$ has order greater or equal to $2$. In particular, keeping
the previous notations, $\Psi_1(\id)\neq\Psi_1(-\id)$.\end{ex}

\bigskip

\subsection{K3 surfaces}\label{subsec:defex}~

\medskip

An \emph{analytic K3 surface} (or simply a \emph{K3 surface}) is a
complex smooth surface $X$ with trivial canonical bundle and such
that $H^1(X,\ZZ)=0$. Note that Siu proved in \cite{Si} that these
properties automatically imply that $X$ carries a K\"{a}hler
structure. Although in this paper we will be mainly interested in
complex non-projective K3 surfaces let us start with some easy
examples of algebraic K3 surfaces.

\begin{ex}\label{ex:K32} {\bf (Complete intersections)} Let $X$ be a smooth surface in $\PP^{n+2}$ which is the complete intersection of $n$
hypersurfaces of degrees $d_1,\ldots,d_n$, with $d_i>1$ for
$i\in\{1,\ldots,n\}$. The adjunction formula shows that the
dualizing sheaf $\omega_X$ of $X$ is $\mathcal{O}_X((d_1+\ldots
+d_n)-(n+3))$. Hence $X$ is a K3 surface if and only if $d_1+\ldots
+d_n=n+3$. In particular the only possibilities are:
\[
\begin{array}{llll}
n=1 & d_1=4; & \\
n=2 & d_1=2 & d_2=3; & \\
n=3 & d_1=2 & d_2=2 & d_3=2.
\end{array}
\]\end{ex}

A result of Kodaira proves that K3 surfaces, considered as
differentiable manifolds, are all diffeomorphic. In particular,
Example \ref{ex:K32} shows that any K3 surface is diffeomorphic to a
smooth quartic surface in $\PP^3$. Since a surface of this type is
simply connected, any K3 surface is simply connected as well.

The \emph{N{\'e}ron-Severi group} of a K3 surface $X$ is the free
abelian group
$$\NS(X):={\rm c}_1(\Pic(X))\subset H^2(X,\ZZ),$$ where ${\rm
c}_1$ is the natural map
\[
H^1(X,\ZZ)\lto H^1(X,\ko_X)\lto \Pic(X)\mor[{\rm c}_1] H^2(X,\ZZ)
\]
induced by the exponential exact sequence
\[
0\lto\ZZ\lto\ko_X\lto\ko_X^*\lto 0.
\]

By the Hodge decomposition theorem, $H^1 (X, \ko_X) = 0$. Therefore ${\rm c}_1$ is injective
and $\Pic(X)\iso\NS(X)$. The rank $\rho(X)$ of the free abelian
group $\NS(X)$ is the \emph{Picard number of} $X$.

The orthogonal complement $T(X):=\NS(X)^\perp\subset H^2(X,\ZZ)$
with respect to the cup product is the \emph{transcendental lattice}
of $X$. Notice that if $X$ is algebraic, the Hodge Index Theorem
implies that $\NS(X)$ has signature $(1,\rho(X)-1)$. In the non-algebraic case this is no longer
true. Indeed, we can consider the following example:

\begin{ex}\label{ex:K33} {\bf (Kummer surfaces)}\index{K3 surface!Kummer surfaces} Let $T$ be a
complex torus of dimension 2 and consider the quotient
$S:=T/\langle\iota\rangle$, where $\iota(a)=-a$ for any $a\in T$.
The surface $S$ has 16 singular points corresponding to the points
fixed by $\iota$. Let $\widetilde T$ be the blow-up of $T$ along
those 16 points and let $\tilde\iota$ be the natural extension of
$\iota$ to $\widetilde T$ (i.e.\ the extension of $\iota$ acting
trivially on the points of the exceptional curves). The quotient
$$\Km(T):=\widetilde T/\langle\tilde\iota\rangle$$ is a K3 surface and
it is called the \emph{Kummer surface of $T$}. Observe
that $\Km(T)$ is the crepant resolution of the singular quotient
$T/\langle\iota\rangle$ and $\Km(T)$ is algebraic if and only if $T$
is an abelian surface. Moreover, if the signature of $\NS(T)$ is
$(s,r)$, then the signature of $\NS(\Km(T))$ is $(s,r+16)$. For
example, the signature $(0,16)$ is realized considering a
$2$-dimensional complex torus $T$ with $\Pic(T)=0$.\end{ex}

The second cohomology group $H^2(X,\ZZ)$ endowed with the cup
product $(-,-)$ is an even unimodular lattice of signature
$(3,19)$ (and hence of rank 22). More precisely, due to Theorem
\ref{thm:milnor}, there exists an isometry
\[
H^2(X,\ZZ)\iso\Lambda:=U^{\oplus 3}\oplus E_8(-1)^{\oplus 2},
\]
where $U$  and $E_8$ are defined in Example \ref{ex:lattices}. The
lattice $\Lambda$ is usually called the \emph{K3 lattice}. Notice that due to this description of $H^2(X,\ZZ)$, $T(X)$ has
signature $(2,20-\rho(X))$ if $X$ is algebraic while it is $(3,19)$ if $X$ is a \emph{generic analytic K3 surface} (i.e.\ $\Pic(X)=0$).

Given a K3 surface $X$, a \emph{marking} for $X$ is the choice of
an isometry $\jj:H^2(X,\ZZ)\isomor\Lambda$. A \emph{marked K3
surface} is a pair $(X,\jj)$, where $X$ is a K3 surface and $\jj$
is a marking for $X$. Two marked K3 surfaces $(X,\jj)$ and
$(X',\jj')$ are isomorphic if there exists an isomorphism
$f:X\isomor X'$ such that $\jj'=\jj\circ f^*$.

\bigskip

\subsection{Hodge structures}\label{subsec:Hodgestr}~

\medskip

Let $X$ be a K3 surface. The usual weight-two Hodge structure on
$H^2(X,\ZZ)$ is the decomposition
\[
H^2(X,\CC)=H^{2,0}(X)\oplus H^{1,1}(X)\oplus H^{0,2}(X).
\]
Since $X$ has trivial dualizing sheaf,
$H^{2,0}(X)=\langle\sigma_X\rangle\iso\CC$, where
\begin{eqnarray}\label{eqn:p}
(\sigma_X,\sigma_X)=0\;\;\;\;\;\mbox{and}\;\;\;\;\;(\sigma_X,\overline{\sigma_X})>0.
\end{eqnarray}
Let $X$ and $Y$ be K3 surfaces. An isometry (with respect to the
cup product) $f:H^2(X,\ZZ)\isomor H^2(Y,\ZZ)$ is a \emph{Hodge
isometry} if $f(H^{2,0}(X))=H^{2,0}(Y)$. An isometry $g:T(X)\isomor T(Y)$ is a Hodge isometry
if $g(H^{2,0}(X))=H^{2,0}(Y)$ (notice that $H^{2,0}(X)\subseteq T(X)\otimes\CC$ and $H^{2,0}(Y)\subseteq T(Y)\otimes\CC$ because
$\NS(X)=H^{1,1}(X)\cap H^2(X,\ZZ)$ and $\NS(Y)=H^{1,1}(Y)\cap
H^2(Y,\ZZ)$).

\begin{remark}\label{rmk:transc} Notice that the Lefschetz
$(1,1)$-theorem implies that the transcendental lattice $T(X)$ of a
K3 surface $X$ can be described as the minimal primitive sublattice
of $H^2(X,\ZZ)$ such that $\sigma_X\in T(X)\otimes\CC$. In particular, if $X$ is generic analytic, then $T(X)=H^2(X,\ZZ)$.\end{remark}

The second cohomology group of a K3 surface contains two interesting
cones. Indeed, the \emph{positive cone} $\kc_{X}\subset
H^{1,1}(X)\cap H^2(X,\RR)$ is the connected component of the cone
$\{x\in H^{1,1}(X)\cap H^2(X,\RR):x\cdot x>0\}$ containing a
K\"{a}hler class. We denote by $\kk_X$ the \emph{K\"{a}hler cone} of
$X$ (i.e.\ the cone generated by the K\"{a}hler classes).

\begin{remark}\label{rmk:pos=kae} One can prove that
$\kk_X=\{x\in\kc_X:(x,\delta)>0:\delta\in\Delta(X)\}$, where
$\Delta(X)=\{y\in\NS(X):(y,y)=-2\}$. It it clear that if $X$ is such
that $\NS(X)=0$, then $\kk_X=\kc_X$.\end{remark}

Of course, the Hodge isometry induced on cohomology by any
isomorphism $f:X\isomor X$ preserves the K\"{a}hler cone. Even more
is true. Indeed the Torelli Theorem shows that the converse of this
statement holds. Namely the Hodge and lattice structures on the
second cohomology group identify a K3 surface up to isomorphism.

\begin{thm}\label{thm:strongtorelli} {\bf (Strong Torelli Theorem)}
Let $X$ and $Y$ be K3 surfaces. Suppose that there exists a Hodge
isometry $g:H^2(X,\ZZ)\isomor H^2(Y,\ZZ)$ which maps a K\"{a}hler
class on $X$ into the K\"{a}hler cone of $Y$. Then there exists a
unique isomorphism $f:X\isomor Y$ such that $f_*=g$.\end{thm}

Many people have given different proofs of the previous result at
different levels of generality. Among them we recall \cite{BR,LP,PS}
(see also \cite{Ge,BPV}). The following is a weaker form of the
previous result:

\begin{thm}\label{cor:torelli} {\bf (Torelli Theorem)} Let $X$ and $Y$ be K3 surfaces. Then
$X$ and $Y$ are isomorphic if and only if there exists a Hodge
isometry $H^2(X,\ZZ)\iso H^2(Y,\ZZ)$.\end{thm}

There is a second way to use the weight-two Hodge structure on
$H^2$ to study the geometry of a K3 surface $X$. Consider
the quasi-projective variety
\[
Q:=\{\sigma\in\PP(\Lambda\otimes\CC):\sigma\cdot\sigma=0\mbox{ and
}\sigma\cdot\overline{\sigma}>0\}.
\]
which is called the \emph{period domain}. If $(X,\jj)$ is a marked
K3 surface, \eqref{eqn:p} implies $[\varphi(\CC\sigma_X)]\in Q$
(here $[\varphi(\CC\sigma_X)]$ is the class in
$\PP(\Lambda\otimes\CC)$ of the line
$\varphi(\CC\sigma_X))\subset\Lambda\otimes\CC$).

Let $\km$ be the moduli space of marked K3 surfaces (see \cite{Ge,Do2} for a discussion about the construction of $\km$). The map
\[
\Theta:\km\lto Q
\]
sending $(X,\varphi)$ to $[\varphi(\CC\sigma_X)]$ is the {\it period
map} and $[\varphi(\CC\sigma_X)]$ is usually called \emph{period} of
the K3 surface $X$.

\begin{thm}\label{thm:surjectivity} {\bf (Surjectivity of the Period Map)} The map $\Theta$ is
surjective.\end{thm}

The interested reader can find a readable account of the proof of this deep result in \cite{BPV}.

Consider a marked K3 surface $(X,\jj)$ such that $\NS(X)\neq0$ and
let $H^{2,0}(X)=\langle\sigma_X\rangle$. In particular, $\rk T(X)<22$ and there exists a
primitive sublattice $L\hookrightarrow\Lambda$ with $\rk L<22$ and
such that $[\varphi(\CC\sigma_X)]\in Q\cap\PP(L\otimes\CC)$. Since
there exist at most countably many primitive sublattices $L$ of
$\Lambda$ with rank strictly smaller than 22, the periods of marked
K3 surfaces $(X,\jj)$ such that $\NS(X)=0$ are contained in the
complement of the union of countably many hyperplane sections of the
period domain $Q$. Therefore, a K3 surface $X$ with $\NS(X)=0$ is
\emph{generic}.

Passing from $H^2(X,\ZZ)$ to the total cohomology group $H^*(X,\ZZ)$
of a K3 surface $X$, we see that $H^*(X,\ZZ)$ has a natural lattice
structure induced by the \emph{Mukai pairing} defined by
\begin{eqnarray}\label{eqn:MukPair}
\langle\alpha,\beta\rangle:=-(\alpha_0,\beta_4)+(\alpha_2,\beta_2)-(\alpha_4,\beta_0),
\end{eqnarray}
for every $\alpha=(\alpha_0,\alpha_2,\alpha_4)$ and
$\beta=(\beta_0,\beta_2,\beta_4)$ in $H^*(X,\mathbb{Z})$. The group
$H^*(X,\mathbb{Z})$ endowed with the Mukai pairing is called
\emph{Mukai lattice} and we write $\widetilde{H}(X,\mathbb{Z})$ for
it. Since $H^0(X,\ZZ)\oplus H^4(X,\ZZ)$ with the cup product is
isometric to $U$, we have a natural isometry of lattices
\[
\widetilde{H}(X,\mathbb{Z})\iso\widetilde\Lambda:=\Lambda\oplus U.
\]

We can now define a weight-two Hodge structure on the total
cohomology group $H^*(X,\ZZ)$ in the following natural way (see
\cite{Hi,Hu,HS1} for a more general setting):
\begin{eqnarray}\label{form:Hodgestru}
\begin{array}{l}
\widetilde H^{2,0}(X):=H^{2,0}(X),\\
\widetilde H^{0,2}(X):=H^{0,2}(X),\\
\widetilde H^{1,1}(X):=H^0(X,\mathbb{C})\oplus H^{1,1}(X)\oplus
H^4(X,\mathbb{C}).
\end{array}
\end{eqnarray}
An isomorphism $f:\widetilde H(X,\ZZ)\isomor\widetilde
H(Y,\ZZ)$ is a Hodge isometry if it is an isometry of lattices which
preserves the Hodge structures just defined.

As in the case of the second cohomology group, for a K3 surface $X$,
the lattice $\widetilde H(X,\ZZ)$ has some very interesting
substructure. Indeed, let $\sigma$ be a generator of $H^{2,0}(X)$
and $\omega$ a K{\"a}hler class (e.g.\ an ample class if $X$ is
algebraic), then
\begin{eqnarray}\index{Positive four-space}
P(X,\sigma,\omega):=\langle{\rm Re}(\sigma), {\rm
Im}(\sigma),1-\omega^2/2,\omega\rangle,
\end{eqnarray}
is a positive four-space in $\widetilde H(X,\RR)$ (here ${\rm
Re}(\sigma)$ and ${\rm Im}(\sigma)$ are the real and imaginary part
of $\sigma$). It comes, by the choice of the basis, with a natural
orientation. It is easy to see that this orientation is independent
of the choice of $\sigma$ and $\omega$.

In general, let $g:\Gamma\to \Gamma'$ be an isometry of lattices
with signature $(4,t)$. Suppose that positive four-spaces
$V\subset\Gamma_\RR:=\Gamma\otimes\RR$ and
$V'\subset\Gamma'_\RR:=\Gamma'\otimes\RR$ and orientations for both
of them have been chosen. Then one says that the isometry $g$
\emph{preserves the given orientation of the positive directions}
(or, simply, that $g$ is \emph{orientation preserving}) if the
composition of
\[
g_{\RR}:V\longrightarrow\Gamma'_\RR
\]
with the orthogonal projection $\Gamma'_\RR\to V'$ is compatible
with the given orientations of $V$ and $V'$. By $\OO_+(\Gamma)$
one denotes the group of all orientation preserving isometries of
$\Gamma$.

A Hodge isometry $f:\widetilde H(X,\ZZ)\isomor\widetilde H(Y,\ZZ)$
is \emph{orientation preserving} if it preserves the orientation of the
four positive directions previously defined.

\begin{ex}\label{ex:orient} (i) It is easy to see that, for any K3 surface $X$, the Hodge
isometry $j:=\id_{H^2(X,\ZZ)}\oplus(-\id)_{H^0(X,\ZZ)\oplus
H^4(X,\ZZ)}$ is not orientation preserving.

(ii) Consider now $\delta=(1,0,1)\in\widetilde H(X,\ZZ)$ and take
the Hodge isometry defined by
\[
s_\delta(x)=x-2\frac{\langle\delta,x\rangle}{\langle\delta,\delta\rangle}\delta=x+\langle\delta,x\rangle\delta,
\]
for any $x\in\widetilde H(X,\ZZ)$. Clearly $s_\delta$ is a
reflection and is orientation preserving.

(iii) Let $f:X\isomor Y$ be an isomorphism of K3 surfaces and take
the natural identifications $H^0(X,\ZZ)=H^0(Y,\ZZ)=\ZZ$ and
$H^4(X,\ZZ)=H^4(Y,\ZZ)=\ZZ$. The Hodge isometry
$f^*|_{H^2(Y,\ZZ)}\oplus(\id)_{\ZZ\oplus\ZZ}:\widetilde
H(Y,\ZZ)\isomor\widetilde H(X,\ZZ)$ preserves the orientation of the
four positive directions.\end{ex}

\bigskip

\subsection{McMullen's examples}\label{subsec:McMul}~

\medskip

In this section we describe some generic analytic K3 surfaces
following McMullen's original idea in \cite{Mc} (see also
\cite[Thm.\ 1.7]{GMc}). This beautiful construction requires some
preliminary results about Salem numbers and Salem polynomials
(mainly contained in \cite{GMc}) which we briefly summarize.

\begin{definition}\label{def:Salem} (i) An irreducible monic polynomial
$S(x)\in\ZZ[x]$ is a \emph{Salem polynomial} if the set of its roots
is
$\{\lambda,\lambda^{-1},a_1,\overline{a_1},\ldots,a_n,\overline{a_n}\}$,
where $\lambda\in\RR$, $\lambda>1$ and $|a_i|=1$, for any
$i\in\{1,\ldots,n\}$.

(ii) A Salem polynomial $S(x)$ is \emph{unramified} if
$|S(-1)|=|S(1)|=1$.

(iii) An algebraic real number $\lambda>1$ is a \emph{Salem number} if the
minimal irreducible monic polynomial $P(x)$ such that
$P(\lambda)=0$ is a Salem polynomial.\end{definition}

By definition a Salem polynomial has even degree. Let $n\geq 3$ be
an odd integer and consider a separable polynomial $C(x)\in\ZZ[x]$
of degree $n-3$ and whose roots are contained in the interval
$[-2,2]$ (recall that $C(x)$ is separable if it does not have
multiple roots in $\CC$). Take now the polynomial
\[
R(x)=C(x)(x^2-4)(x-a)-1,
\]
where $a\in\ZZ$. Finally, consider the polynomial
$S(x)=x^nR(x+x^{-1})$ which has degree $2n$. In \cite{GMc} Gross
and McMullen showed that $S(x)$ is an unramified Salem polynomial for any
$a\gg 0$. Hence the following holds true:

\begin{thm} {\bf (\cite{GMc}, Theorem 1.6.)}\label{thm:infSalPol} For any odd
integer $n\geq 3$ there are infinitely many unramified Salem
polynomials of degree $2n$.\end{thm}

An interesting property of Salem polynomials contained in
\cite[Sect.\ 10]{Mc} and which will be used in Section
\ref{subsec:nonalg} is given by the following lemma.

\begin{lem}\label{lem:finitenumbSalem} Given a positive integer $n$
and two real numbers $c_1<c_2$, there exist finitely many Salem
polynomials $S(x)$ such that $\mathrm{deg}\,S(x)=2n$ and
$c_1<\mathrm{tr}\,S(x)<c_2$, where $\mathrm{tr}\,S(x)$ is the sum of
the roots of $S(x)$.\end{lem}

\begin{proof} By definition, the set of roots of a Salem polynomial
of degree $2n$ is of the form
$$\{\lambda,\lambda^{-1},a_1,\overline{a_1},\ldots,a_{n-1},\overline{a_{n-1}}\},$$
where $\lambda\in\RR$, $\lambda>1$ and $|a_i|=1$, for any
$i\in\{1,\ldots,n-1\}$. In particular,
\[
1-n\leq\sum_{i=1}^{n-1}a_i+\overline{a_i}\leq n-1
\]
and, by assumption,
\[
c_1\leq\lambda+\lambda^{-1}+\sum_{i=1}^{n-1}a_i+\overline{a_i}\leq c_2.
\]
Hence $\lambda$ is bounded from above and all the roots of $S(x)$
are contained in a compact subset of $\CC$.

Let us write $S(x)=x^{2n}+\sum_{k=1}^{2n-1}b_kx^k+1$, where $b_k$ is
an integer which is a symmetric function of the roots of $S(x)$.
Therefore, the coefficients of $S(x)$ are contained in a compact
subset of $\RR$. This concludes the proof.\end{proof}

The main idea in McMullen's construction  is to define a K3
surface $X$ with an automorphism $f\in\Aut(X)$ such that the
minimal polynomial of the induced isometry $f^*:H^2(X,\ZZ)\to
H^2(X,\ZZ)$ is a Salem polynomial.

As a first step in this direction we consider an unramified  Salem
polynomial $S(x)$ of degree $22=\rk\Lambda$ whose existence is
granted by Theorem \ref{thm:infSalPol}. Let $\lambda$ be the unique
root of $S(x)$ such that $\lambda\in\RR$ and $\lambda>1$. Consider
the ring
\[
B:=\ZZ[x]/(S(x))
\]
and let $K$ be its field of fractions.

The ring $B$ contains a subring $R$ generated by the element
$y=x+x^{-1}$ (see \cite[Sect.\ 8]{Mc}). Choosing a unit $u\in R$, we
define on $B$ the following bilinear symmetric form:
\begin{eqnarray}\label{eqn:quadrform}
\langle
f_1,f_2\rangle_{B,u}:=\mathrm{Tr}^K_\QQ\left(\frac{uf_1(x)f_2(x^{-1})}{r'(y)}\right),
\end{eqnarray}
where $r$ is the minimal polynomial of $\lambda+\lambda^{-1}$ and
$r'$ is its derivative. McMullen showed that $\langle
f_1,f_2\rangle_{B,u}\in\ZZ$, for any $f_1,f_2\in B$, and so the pair
$(B,\langle-,-\rangle_{B,u})$ is a lattice.

The choice of the unit $u$ is extremely relevant in McMullen's
construction. In fact, in \cite[Sect.\ 8]{Mc} he proved the
following result:

\begin{prop}\label{prop:McK3} There exists $u\in R$ such that the
lattice $(B,\langle-,-\rangle_{B,u})$ is isometric to the K3 lattice
$\Lambda$.\end{prop}

In particular, the lattice $(B,\langle-,-\rangle_{B,u})$ has the
same lattice structure as the second cohomology group of any K3
surface.

The second key step consists in defining an appropriate weight-two
Hodge structure on $B$. To this end, consider now the isomorphism
$F:B\to B$ such that $F(g)=xg$. By definition $F$ preserves the
bilinear form \eqref{eqn:quadrform} and hence yields an isometry of
the lattice $(B,\langle-,-\rangle_{B,u})$. McMullen observes that
the morphism $F+F^{-1}:B\otimes\CC\to B\otimes\CC$ has a unique
eigenspace $E$ whose signature is $(2,0)$. Moreover $E$ can be
decomposed as the direct sum $E=H_1\oplus H_2$ where $H_1$ is the
eigenspace of an eigenvalue $\delta$ of $F$, while $H_2$ is the
eigenspace of the eigenvalue $\delta^{-1}$. In particular $E$ is
invariant under the action of the isometry $F$. Therefore, the
lattice $(B,\langle-,-\rangle_{B,u})$ has a natural weight-two Hodge
structure defined by
\begin{eqnarray}\label{eqn:HodstrMc}
(B\otimes\CC)^{2,0}:=H_1\;\;\;\;\;(B\otimes\CC)^{0,2}:=H_2\;\;\;\;\;(B\otimes\CC)^{1,1}:=(H_1)^{\perp_{\langle-,-\rangle_{B,u}}}.
\end{eqnarray}
Moreover, if $H_1=\langle\sigma\rangle$, then
$H_2=\langle\overline{\sigma}\rangle$ and
\begin{eqnarray}\label{eqn:HodstrMc2}
\langle\sigma,\sigma\rangle_{B,u}=0\;\;\;\;\;\langle\sigma,\overline{\sigma}\rangle_{B,u}>0.
\end{eqnarray}

Since by Theorem \ref{thm:surjectivity} a K3 surface is determined
(up to isomorphism) by the choice of a weight-two Hodge structure as
in \eqref{eqn:HodstrMc} and \eqref{eqn:HodstrMc2}, there exists a K3
surface $X$ and a Hodge isometry
$(H^2(X,\ZZ),(-,-))\iso(B,\langle-,-\rangle_{B,u})$.

Now we can easily prove (using \cite[Thm.\ 3.4]{Mc}) that the
surface $X$ just constructed is a generic analytic K3 surface.

\begin{lem}\label{lem:Mcgen} Let $X$ be a K3 surface as
above. Then $\NS(X)=0$.\end{lem}

\begin{proof} Due to \cite[Sect.\ 8]{Mc}, the minimal polynomial of the isometry $F:B\isomor B$
is the Salem polynomial $S(x)$. Suppose $\NS(X)\neq0$ and let
$F_\QQ$ be the $\QQ$-linear extension of $F$. Then $F_\QQ=F_1\times
F_2$, where $F_1\in\OO(\NS(X)\otimes\QQ)$ and
$F_2\in\OO(T(X)\otimes\QQ)$. But this is impossible, since $S(x)$ is
irreducible and then $F$ has no rational invariant
subspaces.\end{proof}

By abuse of notation, let us denote by $F\in\OO(H^2(X,\ZZ))$ the
isometry induced by $F:B\to B$.

\begin{lem}\label{lem:isomMc} There exists $f\in\Aut(X)$ such that
$f^*=F:H^2(X,\ZZ)\to H^2(X,\ZZ)$.\end{lem}

\begin{proof} Since $\NS(X)=0$, due to Remark \ref{rmk:pos=kae},
$\kk_X=\kc_X$. Hence, by Theorem \ref{thm:strongtorelli}, it is
enough to show that $F(\kc_X)=\kc_X$.

Let $E(\lambda)$ be the eigenspace of the eigenvalue $\lambda$
(which is a Salem number) and let $\alpha\in\CC$ be such that
$F(\sigma_X)=\alpha\sigma_X$, where
$H^{2,0}(X)=\langle\sigma_X\rangle$. As
$$(\sigma_X,\overline{\sigma_X})=(F(\sigma_X),F(\overline{\sigma_X}))=|\alpha|^2(\sigma_X,\overline{\sigma_X})\neq 0,$$
the eigenvalue $\alpha$ is in the unit circle and $E(\lambda)\subset
H^{1,1}(X)\cap H^2(X,\RR)$. Moreover, since $\lambda$ is real, there
exists $0\neq v\in E(\lambda)\cap H^2(X,\RR)$. As $F$ is an
isometry, $(v,v)=(F(v),F(v))=\lambda^2(v,v)$ and then $(v,v)=0$.
Hence, up to changing $v$ with $-v$ we can suppose $v$ in the
boundary of the positive cone. In particular, $F$ preserves the
closure of the positive cone and, therefore, $\kk_X$.\end{proof}

At the end we produced a K3 surface $X$ with $\NS(X)=0$ and such
that $\Aut(X)$ contains an automorphism $f$ whose induced isometry
$f^*:H^2(X,\ZZ)\isomor H^2(X,\ZZ)$ has characteristic polynomial of
Salem type.

\section{Automorphisms}\label{sec:automorphisms}

In this section we study the automorphism groups of non-algebraic K3
surfaces with Picard number $0$. As it will turn out (Theorem
\ref{thm:sumnonalg}) such a group is either trivial or isomorphic to
$\ZZ$. The latter holds true if and only if the K3 surface is
constructed as in Section \ref{subsec:McMul}. Section
\ref{subsec:alg} is devoted to show a few differences with the
algebraic case.

We will reconsider the automorphisms of a non-projective generic K3
surface under a different light in Section
\ref{sec:autoequivalences}. There we will show the close relation
between the automorphism group and the group of autoequivalences of
the derived category.

\bigskip

\subsection{The non-algebraic generic case}\label{subsec:nonalg}~

\medskip

Since any automorphism $f:X\isomor X$ of a K3 surface $X$ induces a
Hodge isometry $f^*:H^2(X,\ZZ)\isomor H^2(X,\ZZ)$, there exists
$\cc(f)\in\CC^\times$ such that $f^*(\sigma_X)=\cc(f)\sigma_X$,
where $H^{2,0}(X)=\langle\sigma_X\rangle$. Hence we get a natural
map
\begin{eqnarray}\label{eqn:chi}
\chi:\Aut(X)\lto\CC^\times
\end{eqnarray}
sending $f$ to $\cc(f)$.

\begin{lem}\label{lem:rootunit} {\rm (i)} For any $f\in\Aut(X)$ of finite order, the complex number
$\cc(f)$ is a root of unity.

{\rm (ii)} If $\rho(X)=0$, then the morphism $\chi$ is
injective.\end{lem}

\begin{proof} To prove (i), suppose that $f\in\Aut(X)$ has finite order. Then $\cc(\langle
f\rangle)$ is a finite cyclic subgroup of $\CC^\times$. Hence all
the eigenvalues of $f^*$ are roots of unity (see also
\cite{Z}).

For (ii), suppose that there exists $f\in\Aut(X)$ such that
$f^*(\sigma_X)=\sigma_X$. This means that $f^*$ has eigenvalue $1$
and then the primitive sublattice $T:=\{x\in H^2(X,\ZZ):f^*(x)=x\}$
is non-empty. Moreover $\sigma_X\in T\otimes\CC$. By Remark
\ref{rmk:transc}, $T=T(X)=H^2(X,\ZZ)$ and so $f^*=\id$. By Theorem
\ref{thm:strongtorelli}, $f=\id$.\end{proof}

We now prove that the automorphism group of non-algebraic generic K3
surfaces does not contain non-trivial automorphisms of finite order.

\begin{prop}{\bf (\cite{Ni1}, Theorem 3.1.)}\label{prop:nikulin} Let $X$ be a K3 surface
such that $\rho(X)=0$ and let $f\in\Aut(X)$ be of finite order. Then
$f=\id$.\end{prop}

\begin{proof} Due to Lemma \ref{lem:rootunit}, it is enough to prove that if $f$ has finite order then
$\cc(f)=1$. Suppose $\cc(f)\neq 1$ and consider the quotient
$X':=X/\langle f\rangle$. Since $\langle f\rangle$ is finite, $X'$
has a finite number of singular points. Let $g:Y\to X'$ be a
resolution of these singularities. Of course, since $X$ is
K\"{a}hler, $Y$ is K\"{a}hler.

Let $\sigma$ be a holomorphic 2-form on $Y$. We can lift $\sigma$ to
a meromorphic 2-form $\sigma'$ on $X$ via the natural quotient map
$X\to X'$. The form $\sigma'$ is not defined just in a finite number
of points. Hence it can be extended to a holomorphic 2-form
$\sigma''$ on $X$ which, moreover, must be invariant under the
action of $\langle f\rangle$. Since we are supposing $\cc(f)\neq 1$,
$\sigma''$ must be trivial and $\sigma$ must be trivial as well. In
particular $\dim H^{2,0}(Y)=0$ and $Y$ is a smooth algebraic surface
(indeed, $\NS(Y)$ must contain $h$ such that $h^2>0$). This would
imply that $\NS(X)\neq0$ which contradicts our
assumption.\end{proof}

Although the group $\Aut(X)$ cannot contain automorphisms of finite
order if $\rho(X)=0$, there may be $f\in\Aut(X)$ with infinite
order. McMullen's construction in Section \ref{subsec:McMul}
provides examples of this type. Indeed, if $X$ is a K3 surface as in
Section \ref{subsec:McMul}, Lemma \ref{lem:isomMc} shows that there
exists $f\in\Aut(X)$ whose minimal polynomial $P(x)\in\ZZ[x]$ is a
Salem unramified polynomial. Since $P(x)$ has a real eigenvalue
$\lambda>1$, $f$ has infinite order (Proposition
\ref{prop:nikulin}).

\begin{prop} {\bf (\cite{Mc}, Theorem 3.2.)}\label{prop:infautsalem}
Let $X$ be a K3 surface such that $\rho(X)=0$ and let $f\in\Aut(X)$
be of infinite order. Then the minimal polynomial of $f^*$ is a
Salem polynomial.\end{prop}

\begin{proof} First notice that $|\cc(f)|\neq 1$.
Indeed, if $\cc(f)$ is a $n$-th root of unity, then $\cc(f^n)=1$ and
by item (ii) of Lemma \ref{lem:rootunit} $f^n=\id$, which is a
contradiction.

More precisely, $f^*$ has an eigenvalue $\lambda$ such that
$|\lambda|\neq 1$. Indeed, all the eigenvalues of $f^*$ are also
roots of its characteristic polynomial in $\ZZ[x]$. Since
$\det(f^*)=\pm1$, the eigenvalues of $f^*$ with absolute value $1$
would be roots of unity. This contradicts the previous remark.

Let $\lambda$ be an eigenvalue of $f^*$ such that $|\lambda|>1$ (which always exists since there is an eigenvalue $\lambda'$ with $|\lambda'|\neq 1$).
The pull-back $f^*$ preserves the subspace $H^{1,1}(X)\subset
H^2(X,\CC)$. Hence $f^*=F_1\times F_2\in\OO(1,19)\times\OO(2,0)$.
Due to \cite[Lemma 3.1]{Mc}, $h\in\OO(p,q)$ has at most
$\min\{p,q\}$ eigenvalues outside the unit circle. Hence $f^*_\CC$
has at most one eigenvalue $\lambda$ with $|\lambda|>1$. Notice that
the uniqueness also implies that such an eigenvalue is real.

The fact that $\lambda>1$ follows easily from the remark that $f^*$,
being induced by an automorphism, must preserve the cone $\kk_X$
(for more details on this, see the proof of Lemma
\ref{lem:isomMc}).\end{proof}

We are now ready to describe the automorphism group of any K3
surface $X$ with $\rho(X)=0$. In this sense, according to \cite{Og1}, we prove now a weak version of the much stronger result \cite[Thm.\ 1.5]{Og1}.

\begin{thm}\label{thm:autnonalg} {\bf (\cite{Og1}, Theorem 1.5.)} Let
$X$ be a K3 surface such that $\rho(X)=0$. Then either
$\Aut(X)=\{\id\}$ or $\Aut(X)=\langle f\rangle\iso\ZZ$.\end{thm}

\begin{proof} Due to \cite[Prop.\ 4.3]{Og1} the
group $\Aut(X)$ is a finitely generated free abelian group. Now
suppose that $\id\neq f,g\in\Aut(X)$ are two distinct generators.

To shorten the notation we put $F:=f^*$ and $G:=g^*$. Let $S_F(x)$,
$S_G(x)$, $\lambda_F$ and $\lambda_G$ be the corresponding Salem
polynomials and Salem numbers. We already know (Lemma
\ref{lem:rootunit}) that the morphism
\[
\cc:\langle f,g\rangle\lto\CC^\times
\]
is injective. Therefore, since $\CC^\times$ is commutative, $FG=GF$.

For any pair $(m,n)\in\ZZ\times\ZZ$, we denote by $S_{(m,n)}(x)$ the
Salem polynomial of the isometry $F^mG^n$. Without loss of
generality we can suppose that the Salem number of $S_{(m,n)}(x)$ is
$\lambda_F^m\lambda_G^n$ (notice that as a consequence of
Proposition \ref{prop:infautsalem} and of the fact that $FG=GF$, $F$
and $G$ are simultaneously diagonalizable).

There exists two real numbers $1<c_1<c_2$ such that the set
$M_1:=\{(m,n)\in\ZZ\times\ZZ\setminus\{(0,0)\}:c_1<\lambda_F^m\lambda_G^n<c_2\}$
is infinite. Indeed, choose $c_1$ and $c_2$ such that
$\lambda_G<c_2/c_1$ (i.e.\ $\log(\lambda_G)<\log(c_2)-\log(c_1)$).
Hence for $m$ negative there always exists at least a positive
integer $n$ such that
$$\log(c_1)-m\log(\lambda_F)<n\log(\lambda_G)<\log(c_2)-m\log(\lambda_F).$$
This is equivalent to $c_1<\lambda_F^m\lambda_G^n<c_2$.

Due to Lemma \ref{lem:finitenumbSalem}, there exists a Salem
polynomial $S(x)$ such that the set $M_2:=\{(m,n)\in
M_1:S_{(m,n)}(x)=S(x)\}$ is infinite. Since $\cc(f^m\circ g^n)$ is a
root of $S(x)$ for any $(m,n)\in M_2$ and as $S(x)$ has only
finitely many roots, there are $m_1,n_1\in\ZZ$ such that
$\cc(f^{m_1}\circ g^{n_1})=\cc(f^{m_2}\circ g^{n_2})$, for any
$(m_2,n_2)\in M_1$. This means that $\cc(f^{m_1-m_2}\circ
g^{n_1-n_2})=1$ and then $f^{m_1-m_2}\circ g^{n_1-n_2}=\id$ (by
Proposition \ref{prop:nikulin}).

As $\cc(\langle f,g\rangle)$ is a free abelian group of finite rank
(\cite[Prop.\ 4.3]{Og1}) and $f$ and $g$ are generators,
$F^{m_1-m_2}\circ G^{n_1-n_2}$ is the identity if and only if
$m_1=m_2$ and $n_1=n_2$. In particular, $M_2=\{(m_1,n_1)\}$ which is
a contradiction.\end{proof}

The following result shows that the case of K3 surfaces $X$ with
$\NS(X)=0$ and with an automorphism of infinite order is relatively
rare.

\begin{cor}{\bf (\cite{Og1}, Proposition
6.2.)}\label{cor:countnumb} There are at most countably many K3
surfaces $X$ with $\rho(X)=0$ and with $f\in\Aut(X)$ of infinite
order.\end{cor}

\begin{proof} Suppose that $F\in\OO(\Lambda)$ has infinite order and that its minimal polynomial is a Salem
polynomial $S(x)$ of degree $22$. Since $S(x)$ is irreducible, the
eigenvalues of $F$ are distinct and the corresponding eigenspaces
have dimension $1$. Due to Lemma \ref{lem:finitenumbSalem} there are
at most countably many Salem polynomials of degree $22$. Hence there
are only countably many subspaces of dimension $1$ of
$\Lambda\otimes\CC$ which are eigenspaces of the $\CC$-linear
extension of some $F\in\OO(\Lambda)$ of infinite order with minimal
polynomial of Salem type.

To conclude the proof observe that if $f\in\Aut(X)$ has infinite
order, then the characteristic polynomial of $f^*$ is a Salem
polynomial of degree $22$ and that $H^{2,0}(X)$ is an eigenspace of
dimension $1$ of the $\CC$-linear extension of $f^*$.\end{proof}

For the sake of simplicity, we summarize the previous results in the
following compact form:

\begin{thm}\label{thm:sumnonalg} Let $X$ be a K3 surface such that
$\rho(X)=0$. Then either $\Aut(X)=\{\id\}$ or $\Aut(X)\iso\ZZ$.
Moreover the latter happens for countably many K3 surfaces.\end{thm}

\bigskip

\subsection{The algebraic case}\label{subsec:alg}~

\medskip

The main difference between the automorphism groups
of generic non-algebraic and algebraic K3 surfaces is that in the first case
a complete description can be given while in the second one only
partial results are available.
This is because of the richer geometry of algebraic K3
surfaces. Although the many partial results about the automorphism
group of algebraic K3 surfaces are very interesting, the techniques
appearing in their proofs are far beyond the scope of this paper and
so we will not consider them here.

More modestly, the aim of this section is to establish a comparison
between the algebraic and the non-algebraic settings mainly through
a few simple examples.

Given a K3 surface $X$, define the group
\[
\OO_{\mathrm{Hodge}}(X):=\{f\in
\OO(T(X)):f(H^{2,0}(X))=H^{2,0}(X)\}.
\]
The following lemma is essentially \cite[Lemma 4.1]{Og}.

\begin{lem}\label{lem:Hodge} Let $X$ be an algebraic K3 surface with odd Picard
number. Then $\mathrm{O}_{\mathrm{Hodge}}(X)=\{\pm
\mathrm{id}\}$.\end{lem}

\begin{proof} Let $H^{2,0}(X)=\langle\sigma_X\rangle$. We can write $T(X)\otimes\RR=P\oplus N$, where
$P=T(X)\otimes\RR\cap(\CC\sigma_X\oplus \CC\overline{\sigma_X})$,
$N=T(X)\otimes\RR\cap H^{1,1}(X)$ and $\sigma_X$ is a generator of
$H^{2,0}(X)$. By the Hodge Index Theorem, $N$ is negative definite
and any Hodge isometry $h:T(X)\to T(X)$ is such that the
$\CC$-linear extension $h_\CC$ is in $\OO(P)\times\OO(N)$. The
subspaces $P$ and $N$ are definite and thus $h$ is diagonalizable
over $\CC$ with eigenvalues with absolute value 1.

Since $h$ is defined over $\ZZ$, the eigenvalues are algebraic
integers and thus they are roots of unity. In particular, there is a
natural number $n$ such that $h(\sigma_X)=\zeta_n\sigma_X$, where
$\zeta_n$ is a primitive $n$-th root of unity. This implies that the
order of $h$ is $n$. Indeed, the primitive sublattice $T=\{x\in
T(X):h^n(x)=x\}$ is such that $\sigma_X\in T\otimes\CC$. Remark
\ref{rmk:transc} implies $T=T(X)$.

Obviously, $\dim N=\mathrm{rk}T(X)-2$ is odd. Since $h$ restricts to
$N$, there exists a non-trivial $x\in T(X)$ such that either
$h(x)=x$ or $h(x)=-x$. We can show that if $h$ has eigenvalue $1$
then $h=\mathrm{id}$. Indeed, if $h(x)=x$, then
$(\zeta_n\sigma_X,x)=(h(\sigma_X),h(x))=(\sigma_X,x)$. As $n\neq 1$,
$(\sigma_X,x)=0$ and $x\in\NS(X)$. This is a contradiction because
$x\in T(X)$.

Similarly one shows that if there exists an $x\in T(X)$ such that
$h(x)=-x$, then $h=-\mathrm{id}$.\end{proof}

Applying the previous result we get a complete description of the
group of automorphisms of some K3 surfaces with Picard number 1. This
shows a nice analogy with the non-algebraic case analyzed in the
previous section.

\begin{prop}\label{prop:K3genalg} Let $X$ be a projective K3 surface
such that $\NS(X)=\langle\ell\rangle$ and $(\ell,\ell)=2d>2$. Then
$\Aut(X)=\{\id\}$.\end{prop}

\begin{proof} Take $f\in\Aut(X)$. By Lemma \ref{lem:Hodge},
$F:=f^*|_{T(X)}=\pm\id$. On the other hand $G:=f^*|_{\NS(X)}=\id$
because $G$ must preserve $\kk_X$. Since $F$ extends to $f^*$,
Proposition \ref{prop:ext} and Example \ref{ex:latt1} imply $F=\id$
and $f^*=\id$. By the Strong Torelli Theorem (Theorem
\ref{thm:strongtorelli}), $f=\id$.\end{proof}

\begin{ex}\label{ex:PP2} An example of an algebraic K3 surface with Picard number $1$
not satisfying the hypotheses in the previous proposition is
obtained by considering a smooth plane sextic $C\subset\PP^2$. The
double cover of $\PP^2$ ramified along $C$ is an algebraic K3
surface. If $C$ is sufficiently generic, then $\rho(X)=1$. Notice
that $X$ has always an automorphism of order two given by the
involution exchanging the two leaves of the covering.\end{ex}

One can consider other algebraic K3 surfaces $X$ with non-trivial
automorphisms of finite order which do not exist in the generic
non-projective case. For example, there exist K3 surfaces $X$ with
an involution $f\in\Aut(X)$ such that $f^*|_{T(X)}=\id$. Such an
automorphism is usually called \emph{Nikulin involution} (see
\cite{GS,GeS,Mo,Mu5,Ni1}). More generally an automorphism $f\in\Aut(X)$ of finite order and such that $f^*|_{T(X)}=\id$ is called \emph{symplectic}. If a K3 surface $X$ has symplectic automorphisms, then $\rk\NS(X)\geq 9$ (\cite{Ni1}). Just to give some precise references, we recall that the case of K3 surfaces $X$ such that $\Aut(X)$ contains a finite subgroup $G$ whose elements are symplectic automorphisms was first studied in \cite{Ni1} for $G$ abelian and in \cite{Mu5} for $G$ non-abelian.

To produce examples of this type, take $\sigma\in (U\oplus
U)\otimes\CC$ such that
\[
(\sigma,\sigma)=0\;\;\;\;\;\mbox{and}\;\;\;\;\;(\sigma,\overline{\sigma})>0.
\]
Since there is a primitive embedding $U\oplus
U\hookrightarrow\Lambda$, the class of the line $\CC\sigma$ lies in
$Q$ and by the Surjectivity of the Period Map (Theorem
\ref{thm:surjectivity}), there exists a marked K3 surface $(X,\jj)$
such that $\jj(\CC\sigma_X)=\CC\sigma$, where
$H^{2,0}(X)=\langle\sigma_X\rangle$. Notice that, by the choice of
$\sigma$, $T(X)\hookrightarrow U\oplus U$. In particular, there
exists a primitive embedding $E_8(-1)^{\oplus
2}\hookrightarrow\NS(X)$. Theorem 5.7 in \cite{Mo} ensures that $X$
has a Nikulin involution and hence a non-trivial automorphism of
finite order. To be precise the isometry induced on cohomology by this automorphism exchanges the two copies of $E_8(-1)$ and fixes the rest of the K3 lattice. Other examples of K3 surfaces with automorphisms group
isomorphic to $\ZZ/2\ZZ$ are discussed in \cite{GL1}.

Nevertheless, as in the non-algebraic case, the automorphism group
of a K3 surface could be infinite. Indeed, for an algebraic K3
surface $X$ and for any $\delta\in\Delta(X)$ (see Remark
\ref{rmk:pos=kae} for the definition of $\Delta(X)$), define the
reflection in $\OO(H^2(X,\ZZ))$
\[
r_\delta(x):=x+(x,\delta)\delta,
\]
for any $x\in H^2(X,\ZZ)$ (see \cite{Do3} for an extensive use of reflections in the study of automorphism groups of K3 surfaces). Let $W(X)$ be the Weyl group generated by
the reflections $r_\delta$. Since any $\delta\in\Delta(X)$ is
orthogonal to the transcendental lattice $T(X)$, $W(X)$ can be also
seen as a subgroup of $\OO(\NS(X))$.

\begin{thm}\label{thm:PiatShaf} {\bf (\cite{PS})} Let $X$ be an algebraic K3 surface.
Then $\Aut(X)$ is finite if and only if $$\OO(\NS(X))/W(X)$$ is a
finite group.\end{thm}

\begin{ex}\label{ex:infiniteaut} (i) Following \cite{W}, one can prove
that the generic algebraic K3 surface $X$ which is a hypersurface of
bidegree $(1,1)$ in $\PP^2\times\PP^2$ has automorphism group
isomorphic to the infinite dihedral group. Other examples of
projective K3 surfaces with infinite (non-cyclic) automorphisms
group are provided in \cite{GL2} where the case of the double cover
of $\PP^2$ is considered (see Example \ref{ex:PP2}).

(ii) One can also produce K3 surfaces $X$ such that
$\Aut(X)\iso\ZZ$. Indeed, following \cite{Bi}, this is the case if
$X$ is a generic projective K3 surface such that the intersection
form on $\NS(X)$ is represented by a diagonal matrix
\[
\left(\begin{array}{cc} 2nd & 0\\ 0& -2n
\end{array}\right),
\]
where $n\geq 2$ and $d$ is not a square. The proof of this fact
involves a direct calculation of the solutions of some Pell
equations.\end{ex}

\section{Coherent sheaves and derived categories}\label{sec:cohder}

In this section we move towards the second main topic of this paper:
the description of the group of autoequivalences of the derived category of a generic analytic
K3 surface.

As a first step we discuss a few properties of the abelian category of
coherent sheaves on surfaces of this type. Not surprisingly, this
will lead to another key difference with respect to the algebraic
case. Indeed an interesting result of Verbitsky (Theorem
\ref{thm:verb}) shows that two smooth analytic K3 surfaces with
trivial N\'{e}ron-Severi lattices have equivalent abelian categories
of coherent sheaves. Notice that in the algebraic case this does not
happen due to a well-known result of Gabriel (Theorem
\ref{thm:gabriel}). As a consequence of Verbitsky's result, not all equivalences are of
Fourier--Mukai type (Proposition \ref{prop:algnoteqnonalg}). Hence
the celebrated result of Orlov about the structure of equivalences
is no longer true in the non-projective setting.

One of the main topics treated in Section \ref{subsec:coh} is the definition of Chern classes for coherent sheaves on compact complex varieties. These classes will be used in the forthcoming sections.

An extensive introduction to derived categories and derived functors with particular attention to the non-algebraic case is included for the convenience of the reader. In Section
\ref{subsec:Serre} we study compact complex manifolds derived equivalent to a complex K3 surface.

\bigskip

\subsection{Coherent sheaves, Chern classes and locally free resolutions}\label{subsec:coh}~

\medskip

In the algebraic setting the abelian category $\coh(X)$ of a smooth
projective variety is a very strong invariant. Indeed the following
classical result of Gabriel shows that it actually encodes all the
geometric information about the variety.

\begin{thm}\label{thm:gabriel} {\bf (\cite{Ga})} Let $X$
and $Y$ be smooth projective varieties such that
$\coh(X)\iso\coh(Y)$ as abelian categories. Then there exists an isomorphism $X\iso
Y$.\end{thm}

This is not the case for generic analytic K3 surfaces. Indeed a result of Verbitsky shows that coherent sheaves are not very
interesting for varieties of this type. The original proof in
\cite{V4} is quite involved and we plan to simplify it in a future
paper.

\begin{thm}\label{thm:verb} {\bf (\cite{V4}, Theorem 1.1.)} Let $X$
and $Y$ be K3 surfaces such that $\rho(X)=\rho(Y)=0$. Then
$\coh(X)\iso\coh(Y)$ as abelian categories.\end{thm}

The same result was proved in \cite{V5} for generic (non-projective)
complex tori following very closely the approach in \cite{V4}.

Coherent sheaves on non-algebraic K3 surfaces share some useful properties with coherent sheaves on algebraic K3
surfaces. More precisely, due to \cite{Se} any coherent sheaf on a
projective algebraic variety admits a finite locally free
resolution. In \cite{Sc} Schuster proved that this is true for
any compact complex surface. In the special case of K3 surfaces, this reads as follows:

\begin{thm}\label{thm:2} Any coherent sheaf $\ke\in\coh(X)$ on a K3 surface $X$ admits a
locally free resolution of finite length.\end{thm}

Note that this does not hold for any analytic variety as it was shown in \cite{Vo}.
As an easy consequence of Schuster's theorem, it makes perfect sense to define the Chern character of any coherent sheaf $\ke\in\coh(X)$ for a K3 surface $X$. Indeed, one considers a locally free resolution of finite length of $\ke$ and takes the (well-defined) alternating sum of the Chern characters of the locally free sheaves appearing in the resolution. Unfortunately this is not sufficient for the purposes of the forthcoming sections since we will also need to calculate the Chern character of coherent sheaves (or complexes) on the product of two K3 surfaces, where Theorem \ref{thm:2} is not known to hold true. For this we need to define Chern classes in a different way.

\medskip

Let $X$ be, more generally, a compact complex manifold. Chern classes for coherent sheaves on $X$ with values in the Hodge cohomology of $X$ (i.e.\ $\bigoplus_p H^p (X, \Omega^p_X)$) were defined by O'Brian, Toledo, and Tong in \cite{OTT2}, substituting locally free resolutions by what they called \emph{twisted cochain complexes}. They also proved Hirzebruch-Riemann-Roch and Grothendieck-Riemann-Roch formulas for these classes (\cite{OTT1, OTT3, OTT4}). Notice that this definition is better suited to the setting of derived categories, due to the close relationship with Hochschild homology (see e.g.\ \cite{C1, C2}, as pointed out in \cite{Block}).

In our work we will follow a different definition by Green contained in his unpublished Ph.D.-thesis \cite{Green} (for a review of it, see \cite{TT}). In Green's approach these classes are defined via \emph{complexes of simplicial vector bundles} and they take values in the de Rham cohomology of the manifold. Moreover, due to \cite[Thm.\ 2]{TT}, they coincide with Chern classes defined by Atiyah--Hirzebruch by considering locally free resolutions in the real-analytic setting.

We present a few details of the construction but the reader is also strongly encouraged to consult the original works for precise definitions and results. In particular, let $\ke$ be a coherent sheaf of $\ko_X$-modules and consider the sheaf $\ko_X^{\rm an}$ of germs of complex-valued real-analytic functions on $X$. Denote by $\ke^{\rm an}$ the $\ko_X^{\rm an}$-module $\ke \otimes_{\ko_X} \ko_X^{\rm an}$. Then $\ke^{\rm an}$ is a coherent $\ko_X^{\rm an}$-module.

\begin{prop}\label{pro:AHRes}  {\bf (\cite{AH}, Proposition 2.6.)}
Let $\ks$ be a coherent sheaf of $\ko_X^{\rm an}$-modules. Then there exists a resolution of $\ks$ over $\ko_X^{\rm an}$
$$0 \to L_k \to \ldots \to L_0 \to \ks \to 0,$$
where $L_i$ is a locally free $\ko_X^{\rm an}$-module, for all $i$.
\end{prop}

According to \cite{AH}, the previous proposition allows one to define Chern classes of coherent sheaves on $X$ with values in $H^* (X,\QQ)$. All the usual properties of Chern classes hold with this definition. For example: commutativity with flat pull-backs, additivity under exact sequences of sheaves (which in particular implies that Chern classes are actually well-defined at the level of Grothendieck group $K(X):=K(\coh (X))$ of the abelian category $\coh (X)$) and Whitney's formula. Moreover, Riemann-Roch formulas (Theorem \ref{thm:GRR}) for Green's Chern classes and a comparison between the various definitions are contained in \cite{Green} and \cite{TT}. Clearly in the case of either compact complex surfaces or projective complex manifolds all these definitions coincide with the classical one via locally free resolutions in the complex analytic category.

\begin{definition}\label{def:Muvect} For $e\in K(X)$, the \emph{Mukai vector of $e$} is the vector $v(e):=\ch(e)\cdot\sqrt{\td(X)}\in H^*(X,\QQ)$.\end{definition}

Just to consider the easiest example, suppose $E$ to be a locally free sheaf on a K3 surface $X$ with rank $r:=\rk E$, first Chern class $\mathrm{c}_1:=\mathrm{c}_1(E)$ and second Chern class $\mathrm{c}_2:=\mathrm{c}_2(E)$. If $e$ is the class of $E$ in $K(X)$, explicit calculations yields
\[
v(e)=(r,\mathrm{c}_1,\mathrm{c}_1^2/2-\mathrm{c}_2+r).
\]
From this we deduce that if $X$ is a generic analytic K3 surface (i.e.\ $\NS(X)=0$), the Mukai vector of any $e\in K(X)$ is $v(e)=(r,0,s)$. For the rest of this paper, given a smooth compact complex variety $X$ and $\ke\in\coh(X)$, we denote by $[\ke]$ the class of $\ke$ in the Grothendieck group of $X$.

\begin{remark}\label{rmk:qcohsheaves} Since $\coh(X)$ does not
univocally encode the geometry of a compact complex manifold $X$, one should work with larger categories of sheaves.
Unfortunately quasi-coherent sheaves are not the best choice:
many different definitions of quasi-coherent sheaves are present in
the literature and none of them seems to be the definitive one. Moreover the
relations between these definitions have not been completely explored at
the moment (see the discussion at the end of Section 5 of \cite{BBT}
for more details). In Section \ref{subsec:derivcomplexmnflds }
we will work with the abelian category of
$\ko_X$-modules.\end{remark}

\bigskip

\subsection{Derived categories and functors}\label{subsec:derived}~

\medskip

For the convenience of the reader we briefly recall the algebraic
construction of the derived category of an abelian category. Since
we tend to be very brief in this section, the reader interested in
more details and in all the proofs is encouraged to consult
\cite{GM} and \cite{H}. On the other hand the expert reader can
easily skip this rather formal section.

Let $\cat{A}$ be an abelian category. Our main examples will be the
abelian category $\Mod{X}$ of $\ko_X$-modules and its full abelian
subcategory $\coh (X)$ of coherent sheaves, for a complex manifold
$X$. A complex $M^{\bullet}$ over $\cat{A}$ is a pair
$(\{M^i\}_{i\in\ZZ},\{d^i\}_{i\in\ZZ})$ such that $M^i\in\cat{A}$
and $d^i:M^i\to M^{i+1}$ is a morphism such that $d^j\circ
d^{j-1}=0$, for any $j\in\ZZ$. For simplicity, a complex will be
sometimes denoted by $M^\bullet$ or by the pair
$(M^\bullet,d^\bullet)$. It can be also useful to think of it as a
diagram
$$\cdots\mor M^{p-1}\mor[d^{p-1}] M^p\mor[d^p] M^{p+1}\mor\cdots,$$
where $p\in\ZZ$. A morphism of complexes is a set of
arrows $\{f^i\}_{i\in\ZZ}$ as in the following diagram
\[
\xymatrix{\cdots\ar[r]^{d_{M^\bullet}^{i-2}} &
M^{i-1}\ar[d]^{f^{i-1}}\ar[r]^{d_{M^\bullet}^{i-1}} &
M^{i}\ar[d]^{f^{i}}\ar[r]^{d_{M^\bullet}^{i}} &
M^{i+1}\ar[d]^{f^{i+1}}\ar[r]^{d_{M^\bullet}^{i+1}} & \cdots\\
\cdots\ar[r]^{d_{L^\bullet}^{i-2}} &
L^{i-1}\ar[r]^{d_{L^\bullet}^{i-1}} &
L^{i}\ar[r]^{d_{L^\bullet}^{i}} &
L^{i+1}\ar[r]^{d_{L^\bullet}^{i+1}} & \cdots}
\]
and such that, for any $i\in\ZZ$, $f^i\circ
d_{M^\bullet}^{i-1}=d_{L^\bullet}^{i-1}\circ f^{i-1}$. In this way
we obtain the abelian category of complexes $\C(\cat{A})$.

A morphism of complexes is a \emph{quasi-isomorphism} if, for any
$p\in\ZZ$, it induces isomorphisms between all the cohomology
objects
$$H^p (M^\bullet) := \ker d_{M^\bullet}^p / \im d_{M^\bullet}^{p-1}$$ of the complex
$M^{\bullet}$. We define the derived category $\D (\cat{A})$ as the
localization (see \cite{GZ}) of the category of complexes with
respect to the class of quasi-isomorphisms. So we get a category
with the same objects as $\C (\cat{A})$ and in which we formally
invert all quasi-isomorphism.

This quick definition is not very explicit. To get a description of
the morphisms in $\D (\cat{A})$ we should first consider another
category $\K(\cat{A})$ whose construction goes as
follows. First observe that a morphism of complexes $f:
M^{\bullet}\to L^{\bullet}$ is \emph{null-homotopic} if there are
morphisms $r^p:M^p\to L^{p-1}$ such that
$$f^p=d_{L^\bullet}^{p-1}\circ r^p + r^{p+1}\circ d^p_{M^\bullet},$$
for all $p \in \ZZ$. The objects of $\K(\cat{A})$ are the complexes
in $\C(\cat{A})$ while the morphisms are obtained from the ones in
$\C(\cat{A})$ by taking the quotient modulo null-homotopic
morphisms. In $\K (\cat{A})$ the quasi-isomorphisms form a
\emph{localizing class} (in the sense of Verdier) and then it is
possible to have an explicit description of the morphisms in
$\D(\cat{A})$ as ``fractions'', essentially in the same way as for
the product in a localized ring.

The basic properties of derived categories are summarized in the following proposition.

\begin{prop}\label{prop:basicpropdercat} {\rm (i)} The category $\D (\cat{A})$ is additive.

{\rm (ii)} The canonical functor $Q: \C (\cat{A}) \mor \D (\cat{A})$
is additive and makes all quasi-isomorphisms invertible and is
universal among all functors $\C (\cat{A}) \mor \cat{C}$ having this
property, where $\cat{C}$ is any category.

{\rm (iii)} The canonical functor $\cat{A} \mor \D (\cat{A})$, which
maps an object $A \in \cat{A}$ to the complex with all zero besides
$A$ in degree $0$, is fully faithful.
\end{prop}

For the geometric applications we have in mind, we will be
interested in some full subcategories of $\D(\cat{A})$.

\begin{ex}\label{ex:catder}
(i) We denote by $\D^+ (\cat{A})$ (respectively $\D^-(\cat{A})$ and
$\Db(\cat{A})$), the full additive subcategory of complexes
$M^{\bullet}$ such that $M^n=0$, for all $n\gg 0$ (respectively
$n\ll 0$ and $\abs{n}\gg 0$). Such a category is equivalent to the
category of complexes where $H^n(M^\bullet)=0$, for all $n\gg 0$
(respectively $n\ll 0$ and $\abs{n}\gg 0$).

(ii) $\D^*_{\cat{A}'}(\cat{A})$ is the full additive subcategory of
complexes $M^\bullet$ such that $H^n(M^\bullet)\in\cat{A}'$, for all
$n\in\ZZ$, where $\cat{A}'$ is a full \emph{thick} subcategory of
$\cat{A}$ (i.e.\ any extension in $\cat{A}$ of two objects of
$\cat{A}'$ is in $\cat{A}'$) and $*=\emptyset, +,-,\mathrm{b}$. For
example $\coh(X)\subseteq\Mod{X}$ is a thick subcategory.
\end{ex}

Denote by $[1]:\D(\cat{A})\mor\D(\cat{A})$ the \emph{shift functor}
which maps a complex $(M^\bullet,d^\bullet_{M^\bullet})$ to the
complex $(M^\bullet[1],d_{M^\bullet[1]})$, defined by
$(M^\bullet[1])^i = M^{i+1}$ and $d_{M^\bullet[1]} =-d_{M^\bullet}$.
Notice that a morphism of complexes is mapped to the same morphisms
with degrees shifted by $1$. A \emph{standard triangle} in
$\D(\cat{A})$ is a sequence
$$X^\bullet\mor[Q(p)]Y^\bullet\mor[Q(q)]Z^\bullet\mor[\partial \ee]X[1]^\bullet,$$
where $Q:\C(\cat{A})\mor\D(\cat{A})$ is the canonical functor,
$$\ee\ :\ 0\mor X^\bullet\mor[p] Y^\bullet\mor[q] Z^\bullet\mor 0$$
is a short exact sequence of complexes in $\C(\cat{A})$, and
$\partial \ee$ is a special morphism in $\D(\cat{A})$ functorial in
$\ee$ and which lifts the connecting morphism $H^i(Z^\bullet)\mor
H^{i+1}(X^\bullet)$ in the long exact cohomology sequence associated
to $\ee$. A \emph{triangle} in $\D(\cat{A})$ is a sequence
$$X^\bullet\mor Y^\bullet\mor Z^\bullet\mor X^\bullet[1]$$
isomorphic to a standard triangle.

More generally, a \emph{triangulated category} is an additive
category $\cat{T}$ endowed with an autoequivalence $X\mapsto X[1]$
and a class of sequences (called \emph{triangles}) of the form
$$X\mor Y\mor Z\mor X[1]$$ briefly denoted by $(X,Y,Z)$ which are stable under isomorphisms
and satisfy a list of axioms (see, for example, \cite[Ch.\
1]{Neeman} for a nice presentation). If $\cat{T}$ is a triangulated
category, a \emph{full triangulated subcategory} of $\cat{T}$ is a
full subcategory $\cat{U}\subseteq\cat{T}$ such that
$\cat{U}[1]=\cat{U}$ and, whenever we have a triangle $(X,Y,Z)$ of
$\cat{T}$ such that $X$ and $Z$ belong to $\cat{U}$, then there is
$Y'\in\cat{U}$ isomorphic to $Y$. A full triangulated subcategory is
again triangulated in the natural way. For example, all the
subcategories of $\D (\cat{A})$ considered in Example
\ref{ex:catder} are full triangulated subcategories.

To a triangulated category $\cat{T}$ is naturally associated an
abelian group called \emph{Grothendieck group} and denoted by
$K(\cat{T})$. It is defined as the quotient of the free abelian
group generated by the isomorphism classes $[X]$ of objects of
$\cat{T}$ divided by the subgroup generated by the relations
$$[X]-[Y]+[Z],$$
where $(X,Y,Z)$ runs through the triangles of $\cat{T}$.

We have just seen how to pass from an abelian category to a
triangulated one. Later we will need to invert this process. This is
related to the notion of $t$-structure. Indeed, a
\emph{$t$-structure} (see \cite{BD}) on a triangulated category
$\cat{T}$ is a strictly (i.e.\ closed under isomorphism) full
additive subcategory $\cat{F} \subseteq \cat{T}$ such that $\cat{F}
[1] \subseteq \cat{F}$ and, if
$$\cat{F}^{\perp}:=\grf{G\in\cat{T}\,:\,\Hom_{\cat{T}}(F, G) = 0,\,\mbox{ for all }F\in\cat{F}},$$
then for any non-zero $E \in \cat{T}$, there exists a triangle
$$F \mor E \mor G \mor F [1]$$
with $F \in \cat{F}$ and $G \in \cat{F}^{\perp}$. The \emph{heart}
$\cat{A}:=\cat{F}\cap\cat{F}^{\perp}[1]$ of the $t$-structure is a
full abelian subcategory of $\cat{T}$, with the short exact
sequences being precisely the triangles in $\cat{T}$ whose vertices
are objects of $\cat{A}$. A $t$-structure is called \emph{bounded}
if
$$\cat{T} = \bigcup_{i,j} (\cat{F} [i] \cap \cat{F}^{\perp} [j]).$$
Note that for a bounded $t$-structure one has $K(\cat{T}) \cong
K(\cat{A})$.

Let us now briefly recall how one can produce functors
between derived categories starting from functors between the
corresponding abelian categories. We will follows Deligne's approach
as outlined in \cite{Keller}. Let $\cat{A}$ and $\cat{B}$ be abelian
categories and let $F:\cat{A}\mor\cat{B}$ be an additive functor.
Take $M^\bullet\in\D(\cat{A})$ and define a functor
$$({\mathbf{r}}F)(-,M^\bullet):(\D(\cat{B}))^\mathrm{op}\mor{\cat{Mod}\text{-}}\ZZ$$
in the following way: for any $L^\bullet\in\D(\cat{B})$.
\[
({\mathbf{r}}F)(L^\bullet,M^\bullet)=\left\{(f,
s):\begin{array}{l}\mbox{(1)}\;\;f\in\Hom_{\D(\cat{B})}(L^\bullet,F(N^\bullet))\\\mbox{(2)}\;\;
 s\in\Hom_{\K (\cat{A})}(M^\bullet,N^\bullet)\mbox{ is a quasi-isomorphism}\end{array}\right\} / \sim,
\]
where $(f,s)\sim (g,t)$ if and only if there exists $(h, u)$ making
all obvious diagrams commutative. If the functor
$({\mathbf{r}}F)(-,M^{\bullet})$ is representable, then we denote by
${\mathbf{R}}F(M^\bullet)$ the object representing it and we say
that ${\mathbf{R}}F$ is defined at $M^{\bullet}\in\D(\cat{A})$. We
call the functor $\mathbf{R}F$ constructed in this way (and defined
in general only on a subcategory of $\D(\cat{A})$) the \emph{right
derived functor} of $F$. Of course, the above construction can be
applied without changes to arbitrary functors
$F:\K(\cat{A})\mor\K(\cat{B})$. Moreover, one can dually define
\emph{left derived functors}
${\mathbf{L}}F:\D(\cat{A})\mor\D(\cat{B})$ by inverting the
direction of the arrows in the above definition. The link between
the previous definition and the classical one is illustrated by the
following result.

\begin{prop} Suppose that $\cat{A}$ has enough injectives and $M^{\bullet}$
is left bounded. Then ${\mathbf{R}} F$ is defined at $M^{\bullet}$
and we have
$${\mathbf{R}} F(M^{\bullet})= F(I^{\bullet}),$$
where $M^{\bullet} \mor I^{\bullet}$ is a quasi-isomorphism with a
left bounded complex with injective components.\end{prop}

\begin{ex} If $M^{\bullet}=M\in\cat{A}$ is concentrated in degree zero, then
$I^{\bullet}$ may be chosen to be an injective resolution of $M$ and
we find that
$$H^n({\mathbf{R}}F(M))=({\mathbf{R}}^n F)(M),$$
the $n$-th right derived functor of $F$ in the classical sense of
Cartan--Eilenberg (see \cite{Hart2}). For example, if $F=\Hom(N,-)$,
with $N\in\cat{A}$, then, for $M\in\cat{A}$, we have the chain of
natural isomorphisms
\begin{equation}\label{eq:Ext}
\Ext^n_{\cat{A}} (N,M) \cong H^n ({\mathbf{R}} \Hom (N,M)) \cong \Hom_{\D (\cat{A})} (N, M [n]),
\end{equation}
for all $n \in \ZZ$.\end{ex}

An additive functor $F:\cat{T}\mor\cat{T}'$ between two triangulated
categories is called \emph{exact} if there exists a functorial
isomorphism
$$\phi_X : F (X [1]) \mor[\sim] (F(X)) [1]$$
such that the sequence
$$F(X) \mor[F(u)] F(Y) \mor[F(v)] F(Z) \mor[(\phi_X) F(w)] (F(X)) [1]$$
is a triangle of $\cat{T}'$, for any triangle of $\cat{T}$.

For the rest of this paper, we will write $M$ and $f$ instead of
$M^\bullet$ and $f^\bullet$ to denote complexes and morphisms of
complexes. Moreover all functors will be supposed to be exact.

\bigskip

\subsection{Derived categories of complex manifolds}\label{subsec:derivcomplexmnflds }~

\medskip

In this section we study the derived and triangulated categories
associated to the abelian categories of coherent sheaves on complex
manifolds. In particular, $X$ will always be a compact complex
manifold.

\begin{definition}
The \emph{(bounded) derived category} of $X$, denoted by $\Db (X)$
is defined as the triangulated subcategory $\Db_{\coh} (\Mod{X})$ of
$\D (\Mod{X})$ whose objects are bounded complexes of
$\ko_X$-modules with cohomologies in $\coh (X)$.
\end{definition}

\begin{ex}\label{ex:sky} For any closed point $x\in X$, the
skyscraper sheaf $\ko_x$ is in $\Db(X)$. Moreover,
\begin{eqnarray}\label{eqn:morsky}
\Hom(\ko_x,\ko_x[i])\iso\left\{\begin{array}{ll}\CC&\mbox{if
}i\in\{0,2\}\\\CC\oplus\CC&\mbox{if
}i=1\\0&\mbox{otherwise}\end{array}\right.
\end{eqnarray}
Any object $\ke\in\Db(X)$ such that $\Hom(\ke,\ke[i])$ is as in
\eqref{eqn:morsky} is called \emph{semi-rigid}.

Now recall that given an abelian category
$\cat{A}$, $A\in\ob(\cat{A})$ is \emph{minimal} if any non-trivial
surjective morphism $A\to B$ in $\cat{A}$ is an isomorphism. It is
an easy exercise to prove that the minimal objects of $\coh(X)$ are
the skyscraper sheaves.\end{ex}

\begin{ex}\label{ex:spherical} In \cite{ST} a special class of objects was introduced.
An object $\ke^\bullet\in\Db(X)$ is a \emph{spherical object} if
$\ke^\bullet\otimes\omega_X\cong\ke^\bullet$, where $\omega_X$ is
the dualizing sheaf of $X$ and
\[
\Hom(\ke^\bullet,\ke^\bullet[i])\iso\left\{\begin{array}{ll} \CC &
\mbox{if }i\in\{0,\dim X\}\\ 0 & \mbox{otherwise.}\end{array}\right.
\]
In particular if $X$ is a K3 surface, any line bundle (and hence
also $\ko_X$) is a spherical object of $\Db(X)$.\end{ex}

The definition of $\Db (X)$ may seem not appropriate. A better
candidate to work with would be the derived category $\Db (\coh
(X))$.  As we will see in the next
section, unfortunately this category is not the natural place to derive functors and to verify Serre duality. Clearly there is a natural exact functor
\begin{equation}\label{eqn:functcoh}
F : \Db (\coh (X)) \mor \Db (X).
\end{equation}
In particular $\coh (X)$ is the heart of a bounded $t$-structure on $\Db (X)$ and there is a natural identification $K(X) \cong K(\Db (X))$ (note that, via this identification, the Chern classes and the Chern character are well-defined for objects in $\Db (X)$).
The best situation is when this functor is actually an equivalence.
In the algebraic case this is a well-known result (see, for example,
Corollary 2.2.2.1 of Expos\'{e} II in \cite{I}).

\begin{prop}
Let $X$ be a smooth projective complex manifold. Then the functor
$F$ in \eqref{eqn:functcoh} is an equivalence.
\end{prop}

In the purely analytic case the situation is more complicated due to
the absence of a good notion of quasi-coherent sheaves (see Remark
\ref{rmk:qcohsheaves}). But in the surface case the result is still
true via Serre duality (see Section \ref{subsec:Serre}).

\begin{prop}\label{prop:linkcoh} {\bf (\cite{BB}, Proposition 5.2.1.)} Let $X$ be a smooth compact
complex surface. Then the functor $F$ in \eqref{eqn:functcoh}
is an equivalence.\end{prop}

\begin{proof} We denote by $\Ext^n_1$ and $\Ext^n_2$ respectively the Ext-groups
$\Ext^n_{\coh(X)}$ and $\Ext^n_{\Mod{X}}$ in the categories
$\coh(X)$ and $\Mod{X}$, that is, by \eqref{eq:Ext}, the groups
$\Hom (-,-[n])$ in the corresponding derived categories.

The result follows using \cite[Lemma 2.1.3]{I} once we show that for
any $\ke,\kf\in\coh(X)$, the natural morphism
$f_n:\Ext_1^n(\ke,\kf)\to\Ext_2^n(F(\ke),F(\kf))$ is an isomorphism
for any $n$. For $n=0$ this is true by the very definition of the
functor $\Hom$ in both categories. For $n=1$ this property is
trivially verified since the extension of coherent sheaves is still
coherent.

To show that $f_n$ is invertible for $n\neq 1$, we apply once more
\cite[Lemma 2.1.3]{I} which asserts that this happens if for any
$\ke,\kf\in\Mod{X}$ and any $f\in\Ext_2^n(\ke,\kf)$, there
exists a surjective morphism $\ke'\epi\ke$ such that $f$ is sent to
zero by the natural morphism $\Ext_2^n(\ke,\kf)\to\Ext_2^n(\ke',\kf)$.

If $n>2$, Serre duality discussed in Section \ref{subsec:Serre} for
analytic varieties implies that $\Ext_2^n(\ke,\kf)=0$ and so there is nothing to prove. Hence we
just need to consider the case $n=2$.

Let us first show that for any $\ke,\kf\in\coh(X)$ and any $x\in X$,
there is $n$ such that $\Ext^2_2(m_x^n\ke,\kf)=0$, where $m_x$
is the maximal ideal of $x$ in $\ko_X$. Due to Serre duality, we can
just show that for $n\gg 0$ one has $\Hom(\kg,m_x^n\ke)=0$, where
$\kg=\kf\otimes\omega_X$ and $\omega_X$ is the dualizing sheaf of
$X$.

Since $\Hom(\kg,m_x^n\ke)$ is finite dimensional over $\CC$, it is sufficient
to to show that for $a\in\NN$ there is $b>a$ such that
$\Hom(\kg,m^a_x\ke)\neq\Hom(\kg,m^b_x\ke)$. If we pick $f:\kg\to
m_x^a\ke$ then it is an easy check using stalks that there is $b$
such that $\im(f_x)\not\subset m_x^b\ke_x$. Hence
$f\notin\Hom(\kg,m^b_x\ke)$.

Now we take $x,y\in X$ such that $x\neq y$ and we choose $n\in\NN$
such that $\Ext^2_2(m_x^n\ke,\kf)=\Ext^2_2(m_y^n\ke,\kf)=0$. The
canonical map $m_x^n\ke\oplus m_y^n\ke\to\kf$ is clearly surjective.
This concludes the proof.\end{proof}

Going back to the case of generic analytic K3 surfaces $X$, we can now state the following crucial property of $\Db(X)$ which completely classifies its spherical objects. This will be important when we will study stability conditions on the derived category of these K3 surfaces.

\begin{lem}\label{lem:uniquesphK3gen} {\bf (\cite{HMS}, Lemma 3.1.)}
Let $X$ be a generic analytic K3 surface. The trivial line bundle $\ko_X$ is up to shift the only spherical
object in $\Db(X)$.\end{lem}

\bigskip

\subsection{Derived functors on complex manifolds}\label{subsec:derfunc}~

\medskip

Let $X$ and $Y$ be compact complex manifolds. Consider their
corresponding derived categories. As we pointed out in Section
\ref{subsec:derived}, not all functors can be derived. Nevertheless,
the following result shows that all functors with a geometric
meaning are subject to this procedure.

\begin{prop}
Let $f,g:X\to Y$ be morphisms and assume $f$ to be proper. Then the
following functors are defined:
\[
\mathrm{\bf R}\kh
om^\bullet:\Db(X)^{\mathrm{op}}\times\Db(X)\lto\Db(X),
\]
\[
-\stackrel{\bf L}{\otimes}
-:\Db(X)\times\Db(X)\longrightarrow\Db(X),
\]
\[
{\bf L}g^*:\Db(Y)\longrightarrow\Db(X),
\]
\[
{\bf R}f_*:\Db(X)\longrightarrow\Db(Y),
\]
\[
\mathrm{\bf
R}\Hom^\bullet:\Db(X)^{\mathrm{op}}\times\Db(X)\longrightarrow\Db(\cat{vect}),
\]
where we denote by $\cat{vect}$ the abelian category of finite dimensional $\CC$-vector spaces.
\end{prop}
\begin{proof}
According to \cite{Spalt}, the functors listed in the proposition
are all well-defined if we consider the corresponding unbounded
derived categories of modules $\D (\Mod{X})$, $\D (\Mod{Y})$ and $\D
(\cat{Vect})$ (where we denote by $\cat{Vect}$ the abelian category
of $\CC$-vector spaces). To prove the proposition we only have to
show that, by restricting the domain to the corresponding derived
categories of bounded complexes with coherent cohomology, the images
of these functors are still in the bounded derived categories of
complexes with coherent cohomology. But these are well-known results
(see, for example, \cite[Ch.\ 5]{GH}).\end{proof}

To give to the reader a first impression of the relationships
between the derived functors just defined, we list here a few of
their properties whose proofs are, for example, in \cite{Spalt}. As
before, $X$, $Y$ and $Z$ are compact complex manifolds while $f:X\to
Y$ and $g:Y\to Z$ are morphisms. Then there exist the following
natural isomorphisms of functors:
\begin{itemize}
\item Assuming $f$ and $g$ proper, $\mathbf{R}(g_*\circ f_*)\cong\mathbf{R}g_*\circ\mathbf{R}f_*$;
\item $\mathbf{L}(f^*\circ g^*)\cong\mathbf{L}f^*\circ\mathbf{L}g^*$.
\item $\mathbf{R}\Hom^\bullet(\kf,\kg)\cong\mathbf{R}\Gamma(X,\mathbf{R}\kh
om(\kf,\kg))$ as bi-functors.
\item $\mathbf{R}f_*\kh om^\bullet(\kf,\kg)\cong\mathbf{R}\kh
om(\mathbf{R}f_*(\kf),\mathbf{R}f_*(\kg))$.
\item Projection Formula: assuming $f$ proper, $\mathbf{R}f_*(\kf)\stackrel{\bf
L}{\otimes}\kg\cong\mathbf{R}f_*(\kf\stackrel{\bf
L}{\otimes}\mathbf{L}f^*(\kg))$.
\item Assuming $f$ proper, $\Hom^\bullet(\mathbf{L}f^*(\kf),\kg)\cong\Hom(\kf,\mathbf{R}f_*(\kg))$. In other
words, $\mathbf{L}f^*$ is left adjoint to $\mathbf{R}f_*$.
\item $\mathbf{L}f^*(\kf\stackrel{\bf L}{\otimes}\kg)\cong\mathbf{L}f^*(\kf)\stackrel{\bf
L}{\otimes}\mathbf{L}f^*(\kg)$.
\item $\kf\stackrel{\bf L}{\otimes}\kg\cong\kg\stackrel{\bf
L}{\otimes}\kf$ and $\kf\stackrel{\bf L}{\otimes}(\kg\stackrel{\bf
L}{\otimes}\kh)\cong(\kf\stackrel{\bf L}{\otimes}\kg)\stackrel{\bf
L}{\otimes}\kh$.
\item $\mathbf{R}\kh
om^\bullet(\kf,\kg)\stackrel{\bf L}{\otimes}\kh\cong\mathbf{R}\kh
om^\bullet(\kf,\kg\stackrel{\bf L}{\otimes}\kh)$.
\item $\mathbf{R}\kh
om^\bullet(\kf,\mathbf{R}\kh om^\bullet(\kg,\kh))\cong\mathbf{R}\kh
om^\bullet(\kf\stackrel{\bf L}{\otimes}\kg,\kh)$.
\end{itemize}
Moreover, if $u:Y'\rightarrow Y$ is a flat morphism fitting in the
following commutative diagram
\[
\xymatrix{X\times_Y Y'\ar[d]_{g}\ar[r]^{v} & X\ar[d]^{f}\\
Y'\ar[r]^{u} & Y,}
\]
then there exists a functorial isomorphism
\[
u^*\mathbf{R}f_*(\kf)\cong\mathbf{R}g_*v^*(\kf).
\]
Such a property is usually called Flat Base Change Theorem.

\begin{definition}\label{def:twFM} (i) Let $X$ and $Y$
be compact complex manifolds. A functor $F:\Db(X)\to\Db(Y)$ is of
\emph{Fourier--Mukai type} (or a \emph{Fourier--Mukai functor}) if
there exists $\ke\in\Db(X\times Y)$ and an isomorphism of functors
$F\iso\FM{\ke}$, where, denoting by $p:X\times Y\to Y$ and
$q:X\times Y\to X$ the natural projections,
$\FM{\ke}:\Db(X)\to\Db(Y)$ is the exact functor defined by
\begin{eqnarray}\label{eqn:FMT}
\FM{\ke}:=\R p_*(\ke\lotimes q^*(-)).
\end{eqnarray}
The complex $\ke$ is a \emph{kernel} of $F$.

(ii) Two compact complex manifolds $X$ and $Y$ are
\emph{Fourier--Mukai partners} if there
exists a Fourier--Mukai transform
$\Phi:\Db(X)\to\Db(Y)$.\end{definition}

It is a nice exercise to prove that the equivalences in the
following example are of Fourier--Mukai type.

\begin{ex}\label{ex:functors} Here $X$ and $Y$ are always meant to
be compact complex smooth surfaces.

(i) Let $f:X\isomor Y$ be an isomorphism. Then
$f^*:\Db(Y)\isomor\Db(X)$ is an equivalence.

(ii) Let $\ke$ be a spherical object (see Example
\ref{ex:spherical}). Consider the \emph{spherical twist}
$T_\ke:\Db(X)\to\Db(X)$ that sends $\kf\in\Db(X)$ to the cone of
$$\Hom(\ke,\kf)\otimes \ke\to \kf.$$ The kernel of $T_\ke$ is given by the cone of the
natural map $\ke^\vee\boxtimes \ke\to\ko_\Delta$, where
$\ko_\Delta:=\Delta_*\ko_X$ and $\Delta$ is the diagonal embedding
$X\hookrightarrow X\times X$ (notice that here we need the existence
of locally free resolutions as stated in Theorem \ref{thm:2}).

(iii) For any line bundle $L\in\Pic(X)$ we have the \emph{line
bundle twist} $(-)\otimes L:\Db(X)\isomor\Db(X)$.\end{ex}

\begin{thm}\label{thm:GRR}  {\bf (\cite{OTT4, TT})}
Let $X$ and $Y$ be compact complex manifolds and let $f:X\to Y$ be a holomorphic map. Then, for all $\ke \in \Db (X)$, Grothendieck-Riemann-Roch formula holds in $H^* (Y, \QQ)$:
$$\ch (\R f_* \ke) \cdot \td (Y) = f_* (\ch (\ke)\cdot\td (X)).$$
In particular, if $Y$ is a point, for all $\ke \in \Db (X)$, Hirzebruch-Riemann-Roch formula holds:
$$\chi (\ke) := \sum_i (-1)^i \dim H^i (X, \ke) = (\ch (\ke)\cdot\td (X))_{\dim X},$$
where $(-)_{\dim X}$ denotes the $\dim X$-component of a class in $H^* (X, \QQ)$.
\end{thm}

As in the projective case, the Chern class of a coherent sheaf on a K3 surface is integral. The previous theorem can be used to prove that $\ch(\ke)$ (and then $v([\ke])$) is integral for any $\ke\in\Db(X\times X)$, where $X$ is a K3 surface, following line by line the original proof of Mukai in \cite[Lemma 10.6]{H} of the same result in the projective case.
Moreover the Hirzebruch-Riemann-Roch formula can be rewritten for K3 surfaces as
\begin{equation}\label{eq:euler}
\cc(\ke,\kf) :=\sum_i (-1)^i\dim\Hom_{\Db(X)}(\ke,\kf[i]) = -\langle v([\ke]), v([\kf])\rangle,
\end{equation}
for any $\ke,\kf\in\Db(X)$ (here $\langle-,-\rangle$ is the Mukai pairing introduced in \eqref{eqn:MukPair}). The pairing $\chi(-,-)$ will be sometimes called \emph{Euler pairing}.

\medskip

Given two compact complex manifolds $X$ and $Y$, let
$\Eq(\Db(X),\Db(Y))$ be the set of equivalences between the derived
categories $\Db(X)$ and $\Db(Y)$. Moreover we denote by
$\Eq^\mathrm{FM}(\Db(X),\Db(Y))$ the subset of $\Eq(\Db(X),\Db(Y))$
containing all equivalences of Fourier--Mukai type.

The projective case is very special, as stated in the following deep
result due to Orlov (see the main result in \cite{Or} and have a
look to \cite{CS, Ka} for more general statements).

\begin{thm}\label{thm:Orlov} {\bf (\cite{Or})} Let $X$ and $Y$ be smooth projective varieties. Then any
equivalence $\Phi:\Db(X)\isomor\Db(Y)$ is of Fourier--Mukai type and
the kernel is uniquely determined up to isomorphism.\end{thm}

Going back to the case of non-algebraic K3 surfaces, the situation
immediately becomes more complicated. Indeed, reasoning as in the
proof of \cite[Cor.\ 5.23, Cor.\ 5.24]{H} and using Theorem
\ref{thm:verb}, we can prove the following result.

\begin{prop}\label{prop:algnoteqnonalg} Let $X$ and $Y$ be non-isomorphic K3
surfaces such that $\rho(X)=\rho(Y)=0$. Then
$\Eq(\Db(X),\Db(Y))\neq\emptyset$ and
$\Eq^\mathrm{FM}(\Db(X),\Db(Y))\neq\Eq(\Db(X),\Db(Y))$.\end{prop}

\begin{proof} By Theorem \ref{thm:verb}, there exists an equivalence
$\Psi:\coh(X)\isomor\coh(Y)$ which naturally extends by Proposition \ref{prop:linkcoh} to an
equivalence $\Upsilon:\Db(X)\isomor\Db(Y)$. If
$\Eq^\mathrm{FM}(\Db(X),\Db(Y))=\Eq(\Db(X),\Db(Y))$, then $\Upsilon$
is of Fourier--Mukai type and there exists $\ke\in\Db(X\times Y)$ and
an isomorphism of functors $$\Upsilon\iso\FM{\ke}:=\R
p_*(\ke\lotimes q^*(-)).$$

Since the minimal objects of $\coh(X)$ are the skyscraper sheaves
$\ko_x$, where $x$ is a closed point of $X$ (Example \ref{ex:sky}),
$\Upsilon$ maps skyscraper sheaves to skyscraper sheaves. Therefore
$\ke\rest{\{x\}\times Y}$ is isomorphic to a skyscraper sheaf and we
naturally get a morphism $f:X\to Y$ and $L\in\Pic(X)$ such that
$\Upsilon\iso L\otimes f_*(-)$. The morphism $f$ is an
isomorphism since $\Psi$ is an equivalence (\cite[Cor.\ 5.22]{H}).
In particular $X$ and $Y$ would be isomorphic.\end{proof}

Actually in the last section we will show that if $X\not\iso Y$ then
$\Eq^\mathrm{FM}(\Db(X),\Db(Y))=\emptyset$.

As before, if $X$ is a K3 surface such that $\rho(X)=0$, we denote
by $\Aut(\Db(X))$ the group of autoequivalences of $\Db(X)$ and by
$\Aut^\mathrm{FM}(\Db(X))$ the subgroup of $\Aut(\Db(X))$ containing
all Fourier--Mukai autoequivalences.

In the next few pages, to shorten the notations we will often
write a functor and its derived version in the same way (e.g.\ we
write $f_*$ instead of $\mathbf{R}f_*$).

\bigskip

\subsection{Serre functor}\label{subsec:Serre}~

\medskip

For a compact complex manifold $X$, the functor
\begin{eqnarray}
\begin{array}{rcl}
S_{X}:\Db{(X)}&\lto&\Db{(X)}\\
\ke&\longmapsto&\mathcal{F}\otimes\omega_X[\dim (X)],
\end{array}
\end{eqnarray}
where $\omega_X$ is the canonical bundle of $X$, is called
\emph{Serre functor}. The main property of this functor is the
existence of the isomorphism of $\CC$-vector spaces
\[
\Hom_{\Db(X)}(\ke,\kf)\iso\Hom_{\Db(X)}(\kf,S_{X}(\ke))^\vee
\]
which is functorial in $\ke$ and $\kf$ (see \cite[Thm.\ 3.12]{H} for the projective case and
\cite[Prop.\ 5.1.1]{BB} for the general  case).

Since in $\Db(X)$ and $\Db(Y)$, the $\Hom$'s are finite dimensional
complex vector spaces, any equivalence $F:\Db(X)\to\Db(Y)$ (not
necessarily of Fourier--Mukai type) commutes with $S_{X}$ and
$S_{Y}$, i.e.\ there exists an isomorphism $S_{Y}\circ F\iso F\circ
S_{X}$. Indeed, since $F$ is fully faithful,
\[
\Hom(\kf,S_X(\kg))=\mathrm{Hom}(F(\kf),F(S_X(\kg)))\;\;\;\mbox{and}\;\;\;\mathrm{Hom}(\kg,\kf)=\mathrm{Hom}(F(\kg),F(\kf)),
\]
for any $\kf,\kg\in\Db(X)$. Moreover,
\[
\mathrm{Hom}(\kf,S_X(\kg))=\mathrm{Hom}(\kg,\kf)^\vee\;\;\;\mbox{and}\;\;\;\mathrm{Hom}(F(\kg),F(\kf))=\mathrm{Hom}(F(\kf),S_Y(F(\kg)))^\vee.
\]
Hence, we get a functorial (in $\kf$ and $\kg$) isomorphism
\[
\mathrm{Hom}(F(\kf),S_Y(F(\kg))\cong
\mathrm{Hom}(F(\kf),F(S_X(\kg))).
\]
Using the fact that $F$ is an equivalence, we get an isomorphism
between the representable functors $\mathrm{Hom}(-,S_Y(F(\kg)))$ and
$\mathrm{Hom}(-,F(S_X(\kg)))$ and hence an isomorphism
$S_Y(F(\kg)))\cong F(S_X(\kg))$ (see \cite{Br}).

\begin{lem}\label{prop:Serre} Let $X$ be a K3 surface and let $Y$ be a compact complex manifold with an equivalence
$$F:\Db(X)\isomor\Db(Y).$$ Then
$Y$ is a surface with trivial canonical bundle.\end{lem}

\begin{proof} The \emph{order} of the dualizing sheaf $\omega_X$ of a compact complex manifold $X$ is the minimal $k\in\NN\cup\{+\infty\}$
such that $\omega_X^{\otimes k}\cong\mathcal{O}_X$. Following the
proof of \cite[Prop.\ 4.1]{H} and using the property of the Serre
functor previously described, we can show that $\dim X=\dim Y$ and
that the orders of $\omega_X$ and $\omega_Y$ are the same.
\end{proof}

Hence, by the Kodaira classification (see \cite{BPV}), $Y$ is
either a K3 surface or a $2$-dimensional complex torus or a Kodaira surface.

\begin{prop}\label{pro:newproof}
If $Y$ is also assumed to be K\"{a}hler, then $Y$ is a K3 surface with the same Picard number of $X$.
\end{prop}
\begin{proof}
By Example \ref{ex:spherical} we know that $\ko_X$ is a spherical object of
$\Db(X)$. It is enough to prove that a $2$-dimensional complex torus
$T$ does not contain spherical objects.

Following for example \cite[Lemma 15.1]{B}, suppose that $\coh(T)$ contains a spherical sheaf $\ke$ of positive
rank. Since $\ke\dual\otimes\ke=\ko_T\oplus\ke nd_0(\ke)$ (here $\ke
nd_0(\ke)$ is the sheaf of endomorphisms with trace $0$),
$\Hom(\ke,\ke[1])=\Ext^1(\ke,\ke)=H^1(T,\ko_T)\oplus H^1(T,\ke
nd_0(\ke))\neq0$, because $T$ is a torus and
$H^1(T,\ko_T)\iso\CC$. Moreover, any torsion sheaf $\kt$ can be
deformed. Thus $\Ext^1(\kt,\kt)\neq0$ and $\kt$ is not
spherical.

Now let $\ke$ be a spherical object in $\Db(T)$ and let $r$ be the
maximal integer such that $H^r(\ke)\neq 0$. There exists a spectral
sequence
\[
E_2^{p,q}=\bigoplus_{i-j=p}\Ext^q_{\coh(Y)}(H^i(\ke),H^j(\ke))\Longrightarrow\Hom_{\Db(Y)}(\ke,\ke[p+q]),
\]
where $H^j(\ke)$ is the $j$-th cohomology sheaf of $\ke$. For $q=1$,
$p=0$ and $i=j=1$ there is some element surviving in
$\Hom(\ke,\ke[1])$ and so $\ke$ cannot be spherical.

Hence $Y$ is a K3 surface. Let $\kn(X):=K(X)/K(X)^\perp$ and $\kn(Y):=K(Y)/K(Y)^\perp$, where the orthogonal is taken with respect to the Euler pairing \eqref{eq:euler}. The
existence of the Chern character and the Hirzebruch-Riemann-Roch Theorem give the identifications
$\kn(X)=H^0(X,\ZZ)\oplus\NS(X)\oplus H^4(X,\ZZ)$ and
$\kn(Y)=H^0(Y,\ZZ)\oplus\NS(Y)\oplus H^4(Y,\ZZ)$.
Clearly an exact equivalence $\Db (X) \cong \Db (Y)$ induces an isometry (with respect to the Euler pairing) $\kn(X)\iso\kn(Y)$.
\end{proof}

\begin{remark} (i) We will see in the next section (Corollary \ref{cor:FMpart2}) that the previous proposition is general: all Fourier--Mukai partners of a K3 surface $X$ are K3 surfaces, with the same Picard number of $X$. Nevertheless, Proposition \ref{pro:newproof} is interesting in itself. Indeed, it works also for equivalences not of Fourier--Mukai type.

(ii) As an easy consequence, we deduce that $X$ is projective if and only if all K3 surfaces $Y$ with $\Db(X)\cong\Db(Y)$ are projective. Indeed, the signatures of $\kn(X)$ and $\kn(Y)$ are the same. A more precise picture emerges from the recent paper \cite{TV} which we will not discuss here.\end{remark}

\section{Autoequivalences}\label{sec:autoequivalences}

We are now ready to move to the second main problem that we want to
discuss: the description of the group of
autoequivalences of Fourier--Mukai type of the derived category of a generic analytic K3
surface.

The first step will be to study the action induced by any
equivalence of Fourier--Mukai type either on the Grothendieck group
or on the Mukai lattice. Secondly we will use the machinery of
stability conditions introduced by Bridgeland to produce a complete
description of any equivalence.

\bigskip

\subsection{Action on K-theory and cohomology}\label{subsec:actcohom}~

\medskip

In a first step, let $X$ and $Y$ be compact complex manifolds. Let $\Phi_\ke:\Db(X)\rightarrow\Db(Y)$ be a Fourier--Mukai
equivalence with kernel $\ke\in\Db(X\times Y)$. Then the induced isomorphism at the level of Grothendieck groups is given by the morphism $\Phi^K_{[\ke]}:K(X)\rightarrow K(Y)$ defined by
\[
\Phi^K_{[\ke]}(e):=q_*([\ke]\cdot p^*(e)),
\]
where $p:X\times Y\rightarrow X$ and $q:X\times Y\rightarrow Y$ are
the natural projections. In other words, the following diagram commutes
\begin{eqnarray}\label{eqn:KT}
\xymatrix{\Db(X)\ar[d]_{[-]}\ar[rr]^{\Phi_\ke}& & \Db(Y)\ar[d]^{[-]}\\
 K(X)\ar[rr]^{\Phi^K_{[\ke]}}& & K(Y).}
\end{eqnarray}

Using  the definition of Mukai vector, our analysis of Fourier--Mukai equivalences of the bounded derived categories can be further developed. Indeed, for $\ke\in\Db(X)$, one can easily consider the Mukai vector $v([\ke])$ of $\ke$ just by means of Definition \ref{def:Muvect}. Now the morphism $\Phi^K_{[\ke]}:K(X)\rightarrow K(Y)$ gives rise to a map $\Phi^H_{v([\ke])}:H^*(X,\QQ)\rightarrow H^*(Y,\QQ)$ such that
\[
\Phi^H_{v([\ke])}:b\longmapsto q_*(v([\ke])\cdot p^*(b)).
\]
The Grothendieck-Riemann-Roch Theorem shows that the following
diagram commutes:
\begin{eqnarray}\label{eqn:coh}
\xymatrix{K(X)\ar[d]_{v(-)}\ar[rr]^{\Phi^K_{[\ke]}}& &
K(Y)\ar[d]^{v(-)}\\H^*(X,\QQ)\ar[rr]^{\Phi^H_{v([\ke])}}& &
H^*(Y,\QQ).}
\end{eqnarray}

From now on, given a Fourier--Mukai equivalence
$\Phi_\ke:\Db(X)\rightarrow\Db(Y)$, we denote by
\[
\Phi^H_\ke:H^*(X,\QQ)\longrightarrow H^*(Y,\QQ)
\]
the morphism induced on cohomology via the diagrams \eqref{eqn:KT}
and \eqref{eqn:coh} and, to simplify the notation, we write $v(\ke):=v([\ke])$, for any
$\ke\in\Db(X)$.

\begin{lem}\label{lem:FM3} Let $X$ and $Y$ be compact complex manifolds
and let $\Phi_\kf:\Db(X)\rightarrow\Db(Y)$ be a Fourier--Mukai
equivalence whose kernel is $\kf\in\Db(X\times Y)$.  Then
\[
\Phi_\kf^H:H^*(X,\QQ)\lto H^*(Y,\QQ)
\]
is an isomorphism of $\QQ$-vector spaces.\end{lem}

\begin{proof} The proof in the non-algebraic case follows line by line \cite[Prop.\ 5.33]{H}. Since it is quite easy, we reproduce it here for the convenience of the reader. For simplicity, we put $\Phi:=\Phi_\kf$ and $\Phi^H:=\Phi_\kf^H$.

The $\QQ$-linearity is clear form the definition of $\Phi^H$. So it is enough to prove that  $\Phi^H$ is bijective. Since $\Phi$
is an equivalence, its right and left adjoints coincide with the
inverses. Moreover, $\Phi\circ\Phi^{-1}=\id=\Phi_{\ko_\Delta}$,
where $\Delta$ is the image of the diagonal embedding
$i:\Delta\hookrightarrow X\times X$. Indeed, for any $\ke\in\Db(X)$,
{\setlength\arraycolsep{2pt}
\begin{eqnarray*}
\Phi_{\ko_\Delta}(\ke) & =&p_*(\mathcal{O}_\Delta\otimes
q^*(\ke))\nonumber \\
 & =&p_*(i_*\mathcal{O}_X\otimes
q^*(\ke))\nonumber\\
 & =&p_*(i_*(\mathcal{O}_X\otimes i^*q^*(\ke))\nonumber\\
 & =&(p\circ i)_*(q\circ i)^*\ke\\&=&\ke,
\end{eqnarray*}}where $p$ and $q$ are the natural projections, all the functors are
supposed to be derived and the second equality is due to the
Projection Formula (see Section \ref{subsec:derfunc}).

Hence it suffices to prove that $\Phi_{\ko_\Delta}^H=\mathrm{id}$.
Since $\ch(\ko_X)=(1,0,0)$ and the diagonal is isomorphic to $X$,
using the Grothendieck-Riemann-Roch Theorem we get
\[
\ch(\ko_\Delta)\cdot\td(X\times X)=i_*(\td(X)).
\]
Dividing by $\sqrt{\td(X\times X)}$ and using the Projection
Formula, we get
\[
\ch(\ko_\Delta)\cdot\sqrt{\td(X\times X)}=i_*(\td(X)\cdot
i^*((\sqrt{\td(X\times X)})^{-1})).
\]
Since $i^*\td(X\times X)=\td(X)^2$, $i_*(\td(X)\cdot
i^*((\sqrt{\td(X\times X)})^{-1}))=i_*(1,0,0)$ and thus
{\setlength\arraycolsep{2pt}
\begin{eqnarray*}
\Phi^H_{\ko_\Delta}(b) & =&p_*(\ch(\ko_\Delta)\cdot\sqrt{\td(X\times
X)}\cdot q^*(b))\nonumber \\
 & =&p_*(i_*(1)\cdot q^*(b))\nonumber\\
 & =&p_*(i_*(i^*(\cdot
q^*(b))))\\&=&b.
\end{eqnarray*}}The penultimate equality is given by the Projection Formula, while
the last one is trivially due to the fact that $p\circ i=q\circ
i=\mathrm{id}$.\end{proof}

Lemma \ref{lem:FM3} can be made more precise for K3 surfaces.

\begin{thm}\label{thm:coh} Let $X$ and $Y$ be K3 surfaces
and let $\Phi_\kf:\Db(X)\rightarrow\Db(Y)$ be a Fourier--Mukai
equivalence whose kernel is $\kf\in\Db(X\times Y)$. Then $\Phi^H_\kf$ restricts to a Hodge isometry
\[
\Phi_\kf^H:\widetilde H(X,\ZZ)\lto\widetilde H(Y,\ZZ).
\]\end{thm}

\begin{proof} The complete proof of this result (originally due to Orlov, \cite{Or}) is in \cite{H}. Hence we will simply sketch the main steps. As before we put $\Phi:=\Phi_\kf$ and $\Phi^H:=\Phi_\kf^H$.

First notice that, the isomorphism $\Phi^H$ defined over $\QQ$ (Lemma \ref{lem:FM3}) restricts to integral coefficients by the integrality of the Chern character as discussed in Section \ref{subsec:derfunc}. This is an easy exercise which we leave to the reader.

To prove that $\Phi^H_\kf$ preserves the Hodge structure defined in
Section \ref{subsec:Hodgestr}, it is enough to show that
$\Phi^H_\kf(H^{p,q}(X))\subset\bigoplus_{r-s=p-q}H^{r,s}(X)$. The
K\"{u}nneth decomposition of
$c:=\mathrm{ch}(\kf)\cdot\sqrt{\td(X\times Y)}$ is
$\sum\alpha_{u,v}\boxtimes\beta_{r,s}$, where $\alpha_{u,v}\in
H^{u,v}(X)$ and $\beta_{r,s}\in H^{r,s}(Y)$. Since $c$ is algebraic,
the non-trivial terms are the ones such that $u+r=v+s$. If
$\alpha\in H^{p,q}(X)$, then
\[
\Phi^H_\kf(\alpha)=\sum\left(\int_X(\alpha\wedge\alpha_{u,v})\right)\beta_{r,s}\in\bigoplus
H^{r,s}(Y).
\]
In particular the non-zero terms satisfy the condition $u+p=v+q=\dim
X$. Thus $p-q=v-u=r-s$.\end{proof}

This is enough to generalize Proposition \ref{pro:newproof}.

\begin{cor}\label{cor:FMpart2} If $X$ is a K3 surface and $Y$ is a smooth compact complex variety which is a Fourier--Mukai partner of $X$, then $Y$ is a K3 surface and $\rho(X)=\rho(Y)$.\end{cor}

\begin{proof} As a consequence of Lemma \ref{prop:Serre}, we have already proved that $Y$ is either a K3 surface or a $2$-dimensional complex torus or a Kodaira surface. To exclude the last two possibilities observe that if $Y$ is a torus or a Kodaira surface, then  $\dim H^*(X,\QQ)\neq\dim H^*(Y,\QQ)$. But this would contradict Lemma \ref{lem:FM3}.

The Picard numbers of $X$ and $Y$ are the same because, by Theorem \ref{thm:coh}, any equivalence $\Phi:\Db(X)\isomor\Db(Y)$ is such that $\Phi^H(H^{2,0}(X))=\Phi^H(\widetilde H^{2,0}(X))=\widetilde H^{2,0}(Y)=H^{2,0}(Y)$. This yields and isometry $T(X)\iso T(Y)$.\end{proof}

In the generic non-algebraic case the following holds:

\begin{cor}\label{cor:FMpart} Let $X$ and $Y$ be generic analytic K3
surfaces which are Fourier--Mukai partners. Then $X\iso Y$.\end{cor}

\begin{proof} By the previous result an equivalence
$\Phi_\ke:\Db(X)\isomor\Db(Y)$ of Fourier--Mukai type induces a Hodge
isometry $\widetilde H(X,\ZZ)\iso\widetilde H(Y,\ZZ)$. The facts
that the weight-two Hodge structures are preserved and that
$T(X)=H^2(X,\ZZ)$ and $T(Y)=H^2(Y,\ZZ)$ imply that there is a Hodge
isometry $H^2(X,\ZZ)\iso H^2(Y,\ZZ)$. Theorem \ref{cor:torelli}
yields now the desired isomorphism $X\iso Y$.\end{proof}

\begin{remark}\label{rmk:FMpart}
As a consequence, a generic analytic K3 surface does not have
non-isomorphic Fourier--Mukai partners. On the other hand, the
algebraic case is very much different. In fact, in \cite{Og} and
\cite{St} it is proved that for any positive integer $N$, there are
$N$ non-isomorphic algebraic K3 surfaces $X_1,\ldots, X_N$ such that
$\Db(X_i)\iso\Db(X_j)$, for $i,j\in\{1,\ldots,N\}$.\end{remark}

Let us develop further the consequences of Theorem \ref{thm:coh}.
Indeed, consider a (not necessarily projective) K3 surface $X$. If
$\OO_\mathrm{Hodge}(\widetilde H(X,\ZZ))$ denotes the set of Hodge
isometries of the Mukai lattice $\widetilde H(X,\ZZ)$, then Theorem
\ref{thm:coh} yields a group homomorphism
\[
\Pi:\Aut^\mathrm{FM}(\Db(X))\lto\OO_\mathrm{Hodge}(\widetilde
H(X,\ZZ)).
\]
Of course, information about the image of $\Pi$ would give
interesting properties of $\Aut^\mathrm{FM}(\Db(X))$. Notice that a
complete description of the group of autoequivalences in terms of a
group homomorphism analogous to $\Pi$ is available for the closest
relatives of K3 surfaces: abelian varieties (\cite{Or3}).

For a quite long time, the situation for K3 surfaces was widely
unknown and the expectation was summarized by the following
conjecture.

\begin{conjecture}\label{conj:orient} For a K3 surface $X$, the image of $\Pi$
coincides with the subgroup $$\OO_+(\widetilde
H(X,\ZZ))\subseteq\OO_\mathrm{Hodge}(\widetilde H(X,\ZZ))$$ of
orientation preserving Hodge isometries.\end{conjecture}

This conjecture was first stated by Szendr\H{o}i in \cite{Sz}, who
has seen it as a ``mirror'' version of a result of Borcea
(\cite{Borcea}) and Donaldson (\cite{Donaldson}) about the
orientation preserving property for diffeomorphisms of K3
surfaces. Bridgeland (\cite{B}) improved it, by proposing a
description of the kernel of $\Pi$ in terms of the fundamental
group of a domain naturally associated to a manifold parametrizing
stability conditions on the triangulated category $\Db(X)$ (see
Section \ref{subsec:stabcond}). In the algebraic case, some
evidence for the truth of Conjecture \ref{conj:orient} has been
found in \cite{HS1} (see also \cite{HLOY4}). The generic twisted
case is proved in \cite{HMS}. In Section \ref{subsec:equiv} the
conjecture will be proved for generic analytic K3 surface
following the proof in \cite{HMS}. More recently, in \cite{HMS1},
the conjecture was proved for smooth projective (untwisted) K3
surfaces.

\begin{ex}\label{ex:orienisom} (i) The Hodge isometry induced on cohomology by the shift functor
$[1]:\Db(X)\isomor\Db(X)$ is $-\id:\widetilde
H(X,\ZZ)\isomor\widetilde H(X,\ZZ)$ which is orientation preserving.

(ii) Let $f:X\isomor Y$ be an isomorphism of K3 surfaces and
consider the Fourier--Mukai equivalence $f^*:\Db(Y)\isomor\Db(X)$.
Obviously, $(f^*)^H=f^*|_{H^2(X,\ZZ)}\oplus\id_{H^0(X,\ZZ)\oplus
H^4(X,\ZZ)}$ which preserves the natural orientation of the four
positive directions (Example \ref{ex:orient}).

(iii) If $X$ is a K3 surface, the sheaf $\ko_X$ is a spherical object (see Example
\ref{ex:spherical}) and the isometry
induced on cohomology by the spherical twist $T_{\ko_X}$ (Example
\ref{ex:functors}) is the reflection with respect to the vector
$(1,0,1)\in H^*(X,\ZZ)$ described in Example \ref{ex:orient} which
preserves the orientation.\end{ex}

Using the previous examples we can prove the following:

\begin{prop}\label{prop:orient} Let $X$ be a generic analytic K3
surface and let $g\in\OO_+(\widetilde H(X,\ZZ))$. Then there exists
$\Phi\in\Aut^\mathrm{FM}(\Db(X))$ such that $g=\Pi(\Phi)$.\end{prop}

\begin{proof} From the very definition of orientation preserving Hodge
isometries and using the fact that $\rho(X)=\rho(Y)=0$, it follows
easily that $g=h_1\circ h_2$, where $h_1=k_1\oplus\id_{H^2(X,\ZZ)}$
and $h_2=\id_{H^0(X,\ZZ)\oplus H^4(X,\ZZ)}\oplus k_2$, for
$k_1\in\OO(H^0(X,\ZZ)\oplus H^4(X,\ZZ))$ and
$k_2\in\OO(H^2(X,\ZZ))$. Up to composing $g$ with $-\id_{\widetilde
H(X,\ZZ)}$ (which, by Example \ref{ex:orienisom} (i), is induced by
the shift functor) we can also suppose that $h_1$ and $h_2$ are
orientation preserving.

In this situation the Torelli Theorem immediately yields an
automorphism $f:X\isomor X$ such that $h_2=f^*$. On the other hand,
$h_1$ is either the identity or the reflection with respect to the
vector $(1,0,1)$. The second possibility is realized by the spherical
twist (Example \ref{ex:orienisom} (iii)).\end{proof}

\bigskip

\subsection{Stability conditions}\label{subsec:stabcond}~

\medskip

Stability conditions were introduced by Bridgeland in \cite{B2} to
deal with Douglas' notion of $\Pi$-stability (\cite{A,Dou}).
Nevertheless, the relevance of stability conditions and the
topology of the manifolds parametrizing them became immediately
apparent in the geometric context as well. In particular, stability
conditions defined on the bounded derived categories of coherent
sheaves on complex varieties were studied in many interesting
geometric cases (see, for example, \cite{B,B4,Ma,Ok,Ok2}). Very nice
surveys for this are \cite{B1,B5}.

As we mentioned in the previous section, Bridgeland suggested in
\cite{B} that a possible strategy to prove Conjecture
\ref{conj:orient} would involve the topology of the manifold
parametrizing stability conditions on the derived categories of K3
surfaces. In this paper, we will not deal with these topological
problems since in Section \ref{subsec:equiv} we will see that the
study of the equivalences of the derived categories of generic
analytic K3 surface just relies on the notion of stability
condition and does not involve any topological property.

Since, due to Proposition \ref{prop:linkcoh}, there exists a natural
equivalence $\Db(X)\iso\Db(\coh(X))$, throughout this section we
will not distinguish between these two categories. Moreover $X$ will
be a generic analytic K3 surface and $\kn(X)$ will denote the
lattice $H^0(X,\ZZ)\oplus H^4(X,\ZZ)\iso U$. We will restrict our
attention to stability conditions of special type.

\begin{definition}\label{def:stabcond} A \emph{(numerical) stability condition} on $\Db(X)$
is a pair $\sigma=(Z,\kp)$ where
$$\xymatrix{Z:\kn(X)\ar[r]&\CC}$$ is a linear map (the \emph{central charge}) while $\kp(\phi)\subset\Db(X)$ are full
additive subcategories for each $\phi\in\RR$ satisfying the
following conditions:

(a) If $0\ne\ke\in\kp(\phi)$, then $Z(\ke)=m(\ke)\exp(i\pi\phi)$ for
some $m(\ke)\in\RR_{>0}$.

(b) $\kp(\phi+1)=\kp(\phi)[1]$ for all $\phi$.

(c) $\Hom_{\Db(X)}(\ke_1,\ke_2)=0$ for all $\ke_i\in\kp(\phi_i)$
with $\phi_1>\phi_2$.

(d) Any $0\ne\ke\in\Db(X)$ admits a \emph{Harder--Narasimhan
filtration} given by a collection of distinguished triangles
$\ke_{i-1}\to\ke_i\to\ka_i$ with $\ke_0=0$ and $\ke_n=\ke$ such
that $\ka_i\in\kp(\phi_i)$ with
$\phi_1>\ldots>\phi_n$.\end{definition}

Following Bridgeland's terminology, the objects in the category
$\kp(\phi)$ are called \emph{semistable} of phase $\phi$ while the
objects $\ka_i$ in (d) are called the \emph{semistable factors} of
$\ke$. It is a simple exercise to show that they are unique up to
isomorphism. The minimal objects of $\kp(\phi)$ are called
\emph{stable} of phase $\phi$.

Given a stability condition $(Z,\kp)$ and $\ke\in\Db(X)$, we put
$\phi_-(\ke):=\phi_1$ and $\phi_+(\ke):=\phi_n$, where $\phi_1$ and
$\phi_n$ are the phases of the first and the last semistable factor
in the Harder--Narasimhan filtration of $\ke$ (see item (d) of the
previous definition). One defines $\kp((0,1])$ to be the subcategory
of $\Db(X)$ containing the zero object and all $\ke\in\Db(X)$ such
that $0<\phi_-(\ke)\leq\phi_+(\ke)\leq 1$.

\begin{remark}\label{rmk:locallyfinite} (i) All stability conditions studied in this section
are tacitly assumed to be \emph{locally finite} (see \cite{B} for
the precise definition). For our purposes it is enough to know that
this additional condition implies that any semistable object
$\ke\in\kp(\phi)$ admits a finite Jordan--H\"{o}lder filtration,
i.e.\ a finite filtration $0=\ke_0\subset\ldots\subset \ke_n=\ke$ with
stable quotients $\ke_{i+1}/\ke_i\in\kp(\phi)$.

(ii) There is an alternative way to give a stability condition on
$\Db(X)$. Indeed \cite[Prop.\ 5.3]{B2} shows that a stability
condition can equivalently be described by a bounded $t$-structure
on $\Db(X)$ with heart $\cat{A}$ and a stability function
$Z:K(\cat{A})\to\CC$ with Harder--Narasimhan filtrations (recall that a stability function $Z$ is a group homomorphism such that $Z(E)\in\HH\cup\RR_{< 0}$, for all $0\neq E\in\cat{A}$, where $\HH$ is the complex upper
half plane). In particular, given a stability condition $\sigma=(Z,\kp)$,
one can prove that $\cat{A}:=\kp((0,1])$ is the heart (usually
called the \emph{heart of $\sigma$}) of a bounded
$t$-structure.\end{remark}

\begin{definition}\label{def:stabman} The space of all locally finite
numerical stability conditions on $\Db(X)$ is denoted by
$\Stab(\Db(X))$.\end{definition}

In general, it is a highly non-trivial task to show that the space
parametrizing stability conditions on a triangulated category is
non-empty but, once this is achieved, \cite[Thm.\ 1.2]{B2} shows
that such a space is a manifold. Moreover, if the triangulated
category is the bounded derived category of coherent sheaves on some
complex manifold, then the space parametrizing stability conditions
is actually a finite-dimensional complex manifold (\cite[Cor.\
1.3]{B2}).

Hence, let us first show that, in our specific case, $\Stab(\Db(X))$
is non-empty. Consider the open subset
\begin{eqnarray}\label{decompR}
R:=\CC\setminus \RR_{\geq-1}=R_+\cup R_-\cup R_0,
\end{eqnarray}
where $R_+:=R\cap\HH$, $R_-:=R\cap(-\HH)$ and $R_0:=R\cap\RR$ with
$\HH$ denoting the upper half-plane. Given $z=u+iv\in R$, following
\cite{B}, the subcategories $\kf(z),\kt(z)\subset\coh(X)$ are
defined as follows:
\begin{itemize}
\item $\kf(z)$ is the full subcategory of all torsion free sheaves of
degree $\leq v$;
\item$\kt(z)$ is the full subcategory that contains
all torsion sheaves and all torsion free sheaves of degree $>v$.
\end{itemize}
Therefore, if $z\in R_+\cup R_0$, $\kf(z)$ and $\kt(z)$ are
respectively the full subcategories of all torsion free sheaves
respectively torsion sheaves while, if $z\in R_-$, $\kf(z)$ is
trivial and $\kt(z)=\coh(X)$.

Consider the subcategories defined by means of $\kf(z)$
and $\kt(z)$ as follows. If $z\in R_+\cup R_0$, we put
\[
\ka(z):=\left\{\ke\in\Db(X):\begin{array}{l}\bullet\;\;H^i(\ke)=0\mbox{
if }
i\not\in\{-1,0\}\\\bullet\;\;H^{-1}(\ke)\in\kf(z)\\\bullet\;\;H^0(\ke)\in\kt(z)\end{array}\right\}.
\]
On the other hand, let $\ka(z)=\coh(X)$ when $z\in R_-$. Due to
\cite[Lemma 3.3]{B}, $\ka(z)$ is the heart of a bounded
$t$-structure for any $z\in R$.

Now, for any $z=u+iv\in R$ we define the function
\[
\begin{array}{rcl}
Z:\ka(z)&\lto&\CC\\\ke&\longmapsto&\langle
v(\ke),(1,0,z)\rangle=-u\cdot r-s-i(r\cdot v),
\end{array}
\]
where $v(\ke)=(r,0,s)=(r,0,r-{\rm c}_2(\ke))$ is the Mukai vector of
$\ke$.

The main properties of the pair $(\ka(z),Z)$ just defined are
summarized by the following proposition.

\begin{prop}\label{lem:stabcond} {\bf (\cite{HMS}, Lemma 3.2 and Proposition 3.4.)}
For any $z\in R$ the function $Z$ defines a stability function on
$\ka(z)$ which has the Harder--Narasimhan property.
\end{prop}

Using item (ii) of Remark \ref{rmk:locallyfinite}, we conclude that
for any $z\in R$ we have a stability condition $\sigma_{z}$ obtained
from the $t$-structure associated to the pair $(\kf(z),\kt(z))$ and
the stability function $Z$ on $\ka(z)$. In \cite[Prop.\ 3.4]{HMS} it
is also proved that such a stability condition is locally finite.
Hence $\sigma_z\in\Stab(\Db(X))$ and $\Stab(\Db(X))$ is non-empty.
Some additional interesting topological properties of
$\Stab(\Db(X))$ (not used in the rest of this paper) are determined
by the following result.

\begin{thm}\label{thm:simplicon} {\bf (\cite{HMS}, Theorem 3.)} The space
$\Stab(X)$ is connected and simply-connected.\end{thm}

\begin{remark}\label{rmk:actaut} As
pointed out in \cite{B2}, there is a left action of the group
$\Aut^\mathrm{FM}(X)$ on the manifold $\Stab(\Db(X))$. Indeed,
given $\sigma=(Z,\kp)\in\Stab(\Db(X))$ and
$\Phi\in\Aut^\mathrm{FM}(X)$ we define $\Phi(\sigma)$ to be the
stability condition $(Z',\kp')\in\Stab(\Db(X))$ where
$Z'=Z\circ(\Phi^H)^{-1}$ and
$\kp'(\phi)=\Phi(\kp(\phi))$.\end{remark}

Using Lemma \ref{lem:uniquesphK3gen} one can show the following remarkable property of $\Stab(\Db(X))$ which will be used in Section \ref{subsec:equiv}.

\begin{prop}\label{prop:sphact} {\bf (\cite{HMS}, Proposition 1.18 and Corollary 1.19.)}
For any $\sigma\in\Stab(\Db(X))$, there is $n\in\ZZ$ such that
$T^n_{\ko_X}(\ko_x)$ is stable in $\sigma$, for any
closed point $x\in X$.\end{prop}

\bigskip

\subsection{Equivalences}\label{subsec:equiv}~

\medskip

For generic analytic K3 surfaces $X$, we obtain a complete
description of the group of autoequivalences of $\Db(X)$ using the
notion of stability condition. The proof immediately shows that the
topology of $\Stab(X)$ is irrelevant in this case.

\begin{prop}\label{thm:autogen}{\bf (\cite{HMS}, Lemma 3.9.)} Let $X$
and $Y$ be complex K3 surfaces such that $\rho(X)=\rho(Y)=0$. If
$\Phi_\ke:\Db(X)\isomor\Db(Y)$ is an equivalence of Fourier--Mukai
type, then  up to shift
$$\Phi_\ke\iso T_{\ko_Y}^n\circ f_*$$
for some $n\in\ZZ$. Here, $T_{\ko_Y}$ is the spherical twist with
respect to $\ko_Y$ and $f:X\isomor Y$ is an isomorphism.\end{prop}

\begin{proof}Suppose $\sigma=\sigma_{(u,v=0)}$ is one of the distinguished
stability conditions constructed in Section \ref{subsec:stabcond}
with $(u,v)\in R$ and $v=0$. Let $\tilde\sigma:=\Phi_\ke(\sigma)$.

Denote by $T$ the spherical twist $T_{\ko_Y}$. By Proposition
\ref{prop:sphact}, there exists an integer $n$ such that all
skyscraper sheaves $\ko_x$ are stable of the same phase in the
stability condition $T^n(\tilde\sigma)$.

The composition $\Psi:=T^n\circ\Phi_\ke$, is again an equivalence of
Fourier--Mukai type, i.e.\ $\Psi:=\Phi_\kf$, for some kernel
$\kf\in\Db(X\times Y)$. Moreover it sends the stability condition
$\sigma$ to a stability condition $\sigma'$ for which all skyscraper
sheaves are stable of the same phase.

Up to shifting the kernel $\kf$ sufficiently, we can assume that
$\phi_{\sigma'}(\ko_y)\in(0,1]$ for all closed points $y\in Y$.
Thus, the heart $\kp((0,1])$ of the $t$-structure associated to
$\sigma'$ (see Remark \ref{rmk:locallyfinite}), which under
$\Phi_\kf$ is identified with $\ka(u)$, contains as stable objects
the images $\Psi(\ko_x)$ of all closed points $x\in X$ and all point
sheaves $\ko_y$.

In \cite[Rmk.\ 3.3]{HMS} it was observed that the only semi-rigid
stable objects in $\ka(u)$ are the skyscraper sheaves. Hence, for
all $x\in X$ there exists a point $y\in Y$ such that
$\Psi(\ko_x)\iso\ko_y$.

This suffices to conclude that the Fourier--Mukai equivalence
$\Psi_\kf$ is a composition of $f_*$ for some isomorphism
$f:X\isomor Y$ and a line bundle twist (see \cite[Cor.\ 5.23]{H}). But there are no non-trivial line bundles on $Y$.\end{proof}

This immediately yield the following complete description of all
Fourier--Mukai autoequivalences in the generic case.

\begin{thm}\label{prop:autgennonproj} {\bf (\cite{HMS}, Theorem 4.)} If $X$ is a K3 surface with $\Pic(X)=0$, then
$$\Aut^\mathrm{FM}(\Db(X))\iso\ZZ\oplus\ZZ\oplus\Aut(X).$$ The first two
factors are generated respectively by the shift functor and  the
spherical twist $T_{\ko_X}$.\end{thm}

This nice description shows that the interesting (and possibly
non-trivial) part of $\Aut^\mathrm{FM}(\Db(X))$ is given by the
group of automorphisms of $X$. Due to Theorem \ref{thm:sumnonalg},
such a group is usually trivial while it is isomorphic to $\ZZ$ for
countably many generic K3 surfaces characterized by the fact that
$\Aut(X)$ is generated by an automorphism $f$ such that the induced
isometry $f^*$ on the second cohomology group has minimal polynomial
of Salem type. Hence, in the complex analytic setting,
$\Aut^\mathrm{FM}(\Db(X))$ is either isomorphic to $\ZZ^{\oplus 2}$
(for K3 surfaces not of McMullen type) or isomorphic to
$\ZZ^{\oplus 3}$ (for the countably many K3 surfaces of McMullen
type).

We conclude observing that, together with Proposition
\ref{prop:orient}, Theorem \ref{prop:autgennonproj} proves
Conjecture \ref{conj:orient} for analytic generic K3 surfaces. More precisely, there exists a short exact sequence of groups
\[
0\lto\ZZ\oplus\ZZ\lto\Aut^\mathrm{FM}(\Db(X))\lto\OO_+(\widetilde H(X,\ZZ))\lto 0,
\]
where the two copies of $\ZZ$ are generated by the shift by $2$ and by $T_{\ko_Y}^2$ respectively.

\begin{remark} The group $\Aut^\mathrm{FM}(\Db(X))$ encodes some interesting geometric features of the K3 surface $X$. Indeed, if $\Pic(X)\neq 0$, this group is larger than the corresponding one for $X$ generic analytic. Indeed, as we have just seen, the equivalences coming from tensoring with non-trivial line bundles are not present in the latter case.

We have already observed that generic analytic K3 surfaces of McMullen type can be distinguished from the ones not of McMullen type using the previous description of the group of autoequivalences.\end{remark}

\medskip

{\small\noindent{\bf Acknowledgements.} The project was completed when the second named author was visiting the Institute for Mathematical Sciences of the Imperial College (London), whose hospitality is gratefully acknowledged. The stay at the Institute was funded by a grant of the Istituto Nazionale di Alta Matematica (Italy). The results about the equivalences of the derived categories of generic analytic K3 surfaces presented in this paper were obtained by the authors in collaboration with Daniel Huybrechts. It is a great pleasure to thank him for the fruitful collaboration. Thanks also to Alice Garbagnati, Bert van Geemen, and Sukhendu Mehrotra who read a preliminary version of this paper.}

\end{document}